\documentclass[%
reprint,
superscriptaddress,
footinbib,
amsmath,amssymb,
aps,
pre,
showkeys,
10pt
]{revtex4-1}

\usepackage{amsfonts, amsbsy, graphicx, subfigure}
\usepackage{dcolumn}% Align table columns on decimal point
\usepackage{bm}% bold math

\usepackage[hidelinks]{hyperref}
\hypersetup{
	colorlinks=true,
	linkcolor=black,
	filecolor=cyan,      
	urlcolor=blue,
}

\begin{document}

% \preprint{APS/123-QED}

\title{Finding NHIM in 2 and 3 degrees-of-freedom with H{\'e}non-Heiles type potential}
%\title{Finding NHIM in a coupled harmonic potential in 2 and 3 degrees-of-freedom}% Force line breaks with \\

\author{Shibabrat Naik}
\email{s.naik@bristol.ac.uk}
\author{Stephen Wiggins}%
\affiliation{%
School of Mathematics, University of Bristol \\
University Walk, Clifton BS8 1TW, Bristol, UK
}%

\begin{abstract}
We present the capability of Lagrangian descriptors for revealing the high 
dimensional phase space structures that are of interest in nonlinear Hamiltonian systems with index-1 saddle. These phase space structures include normally hyperbolic invariant manifolds (NHIM) and their stable and unstable manifolds, and act as codimenision-1 barriers to phase space transport. The method is applied to classical two and three degrees-of-freedom Hamiltonian systems which have implications for myriad applications in physics and chemistry.
\end{abstract}

\maketitle

% \tableofcontents

\section{\label{sec:intro}Introduction}

It is well-known now that the paradigm of escape from a potential well and the topology of phase 
space structures that mediate such escape 
are used in a broad array of problems such as isomerization of molecular 
clusters~\cite{Komatsuzaki2001}, reaction rates  in chemical 
physics~\cite{Komatsuzaki1999,WiWiJaUz2001}, ionization of a hydrogen atom under electromagnetic 
field in atomic 
physics~\cite{JaFaUz2000}, transport of defects in solid state and semiconductor 
physics~\cite{Eckhardt1995}, buckling modes in structural 
mechanics~\cite{Collins2012,ZhViRo2018}, ship motion and 
capsize~\cite{Virgin1989,ThDe1996,NaRo2017}, escape and recapture of comets and 
asteroids in celestial mechanics~\cite{JaRoLoMaFaUz2002,DeJuLoMaPaPrRoTh2005,Ross2003}, and escape 
into inflation or re-collapse to 
singularity in cosmology~\cite{DeOliveira2002}. As such a method that can identify the high 
dimensional phase space 
structures using low dimensional surface as probes can aid in quantifying the 
escape rates. These low 
dimensional surfaces has been shown to be of as {\it reactive islands} in chemical physics and lead 
to insights into 
sampling rare transition events~\cite{patra_classical-quantum_2015,patra_detecting_2018}. However, 
to benchmark the 
methodology, we first applied it to linear systems where the closed-form analytical expression of 
the phase space 
structures is known~\cite{naik2019finding}. As the next step, in this article, we will focus on 
nonlinear Hamiltonian 
systems which have been extensively studied as ``built by hand'' models of galactic dynamics and for demonstrating quantum dynamical tunneling~\cite{barbanis_isolating_1966,brumer_variational_1976,davis_semiclassical_1979,heller_molecular_1980,waite_mode_1981,kosloff_dynamical_1981,contopoulos_simple_1985,founargiotakis_periodic_1989,barbanis_escape_1990,babyuk_hydrodynamic_2003}. The nonlinear Hamiltonian systems considered here have an underlying H{\'e}non-Heiles type 
potential with the simplest form of nonlinearity, and show regular, quasi-periodic, and chaotic trajectories along with 
bifurcations of periodic orbits. A H{\'e}non-Heiles type 
potential has a well with bottlenecks connecting the region of bounded motion (trapped region) to unbounded motion (escape off to infinity), and have rotational symmetry. 
In addition, these H{\'e}non-Heiles type 
potentials are studied as first benchmark nonlinear systems in applying new phase space transport methods to astrophysical and molecular motion. 
% Another aspect of the coupled harmonic potentials with bottlenecks to infinity is that chaotic trajectories can show both trapping and scattering behavior which are typical of rare escape trajectories. 
% The role of phase space structures in sampling these trajectories has been explained in the context of closed potential wells for two 
% degrees-of-freedom system ~\cite{patra_detecting_2018,patra_classical-quantum_2015}. 
In this article, we will present verification of a method that uses trajectory diagnostic on a low dimensional surface for revealing the phase 
space structures in 4 or more dimensions.

Conservative dynamics on an open 
potential well has received 
considerable attention because the phase space structures, normally hyperbolic invariant manifolds (NHIM) and its invariant 
manifolds, explain the intricate 
fractal structure of ionization 
rates~\cite{mitchell_geometry_2003_I,mitchell_geometry_2003_II,mitchell_chaos-induced_2004}. 
Furthermore, the discrepancies in observed and predicted ionization rates in atomic systems has 
also been explained by 
accounting for the topology of the phase space structures. These have been connected with the 
breakdown of ergodic 
assumption that is the basis for using ionization and dissociation rate 
formulae~\cite{de_leon_intramolecular_1981}. This 
rich literature on chaotic escape of electrons from atoms sets a precedent for applying new methods 
for finding NHIM and its 
invariant manifolds in Hamiltonian with open potential wells 
~\cite{mitchell_analysis_2004,mitchell_chaos-induced_2004,mitchell_nonlinear_2009,mitchell_structure_2007,wang_photoionization_2010}.
 
As we noted earlier, trajectory diagnostic methods which can probe phase space to detect the high dimensional invariant manifolds have potential to be of use in many degrees-of-freedom models. One such method is the Lagrangian descriptors (LDs) that can reveal phase space structures 
by encoding geometric property of trajectories (such as, phase space arc length, configuration space distance or 
displacement, cumulative action 
or kinetic energy) initialised on a two dimensional 
surface~\cite{madrid2009,mendoza2010,mancho2013,lopesino2017}. 
The method was originally developed in the context of Lagrangian transport in time-dependent two 
dimensional fluid 
mechanics. However, it has also been successful in locating transition state trajectory in chemical 
reactions~\cite{balibrea2016lagrangian,craven2017lagrangian,junginger2016lagrangian}. Besides, 
also being applicable to both Hamiltonian and non-Hamiltonian systems, as well as to systems with  arbitrary time-dependence such as stochastic and dissipative forces, and geophysical data from satellite and numerical simulations~\cite{amism11,mendoza2014,ggmwm15,lopesino2017,ramos2018}.

The method of Lagrangian descriptor (LD) is straightforward to implement computationally and it 
provides a ``high resolution'' method for exploring the influence of high dimensional phase space structure on 
trajectory behaviour.  The 
method of LD takes an \emph{opposite} approach to that of classical Lyapunov exponent type 
calculations by 
emphasizing the initial conditions of trajectories, rather than their advected locations that is 
involved in calculating 
normalized rate of divergence. This is achieved by considering a two dimensional section of the 
full phase space and discretizing with a dense 
grid of initial conditions. Even though the trajectories wander off in the phase space, as the 
initial conditions evolve in 
time, there is no loss in resolution of the two dimensional section. In contrast to inferring the 
phase space structures from Poincar\'e sections, LD plots do not suffer from loss of resolution 
since the affects of the structure are encoded in the initial conditions and there is no need for 
the trajectory to return to the section. Our objective is to clarify the use of Lagrangian 
descriptors as a diagnostic on two dimensional sections of high dimensional phase space structures. 
This diagnostic is also 
meant to be used as the preliminary step in computing the NHIM, their stable and unstable manifolds 
using other 
computational means~\cite{junginger2016transition,bardakcioglu2018,ezra_2018}. In this article, we 
will present the method's 
capability to detect the high dimensional phase space structures such as the NHIM, their stable, 
and unstable manifolds in 2 
and 3 DoF Hamiltonian systems. 

\section{Models and Method}

\subsection{\label{sec:model_prob_2dof}Model system: coupled harmonic 2 DoF Hamiltonian}

As pointed out in the Introduction, our focus is to adopt a well-understood model system which is a 2 degrees-of-freedom coupled harmonic oscillator with the Hamiltonian
\begin{equation}
\begin{aligned}
\mathcal{H}(x,y,p_x,p_y) =& T(p_x, p_y) + V_{\rm B}(x,y) \\ 
=& \frac{1}{2}p_x^2 + \frac{1}{2}p_y^2 + \frac{1}{2}\omega_x^2 x^2 + \frac{1}{2}\omega_y^2 y^2 +
\delta x y^2  
\end{aligned}
\label{eqn:Hamiltonian_Barbanis}
\end{equation}
%
%potential energy
%and the model Hamiltonian is given by
%% 
%\begin{equation}
%  V_{\rm B}(x,y) = \frac{1}{2}\omega_x^2 x^2 + \frac{1}{2}\omega_y^2 y^2 + \delta x y^2  
%\label{eqn:Hamiltonian_Barbanis}
%\end{equation}
%
where $\omega_x, \omega_y, \delta$ are the harmonic oscillation frequencies of the $x$ and $y$ 
degree-of-freedom, and the coupling strength, respectively. We will fix the parameters as $\omega_x 
= 1.0, \omega_y = 1.1, \delta = -0.11$ in this study. The two degrees-of-freedom potential is 
also referred to as {\it Barbanis} potential, and has been investigated as a model of galactic 
motion~(\cite{contopoulos1970,barbanis_isolating_1966}), dynamical tunneling and molecular spectra 
in physical chemistry~(\cite{heller1980,davis1981, martens1987}), structural mechanics and ship 
capsize~(\cite{ThDe1996,NaRo2017}).

The equilibria of the Hamiltonian vector field are located at 
\begin{equation}
\left(-\frac{\omega_y^2}{2\delta}, \pm 
\frac{1}{\sqrt{2}}\frac{\omega_x \omega_y}{\delta}, 0, 0 \right) \qquad \text{and}  \qquad \left(0, 
0, 0, 0 \right)
\end{equation}
and are at total energy $E_c = \frac{\omega_x^2 \omega_y^4}{8 \delta^2}$ and $0$ respectively.  The 
energy of the two index-1 saddles (as defined and shown in App.~\ref{ssect:linear})  located 
at 
positive and 
negative y-coordinates and positive 
x-coordinate for $\delta < 0$ will be referred to as {\it critical energy}, $E_c$. In our 
discussion, we will refer to the total energy of a trajectory or initial condition in terms 
of the excess energy, $\Delta E = E_c - e$, which can be negative or positive to denote energy 
below or above the critical energy. For the parameters used in this study, the index-1 saddle 
equilibrium points are located at $\left( 5.5, \pm 7.071, 0, 
0 \right)$ and have energy, $E_c = 15.125$.
%$\left( 10.28992, \pm 5.29412, 0, 0 \right)$

\begin{figure*}[!th]
	\centering
	\subfigure[]{\includegraphics[width=0.3\textwidth]{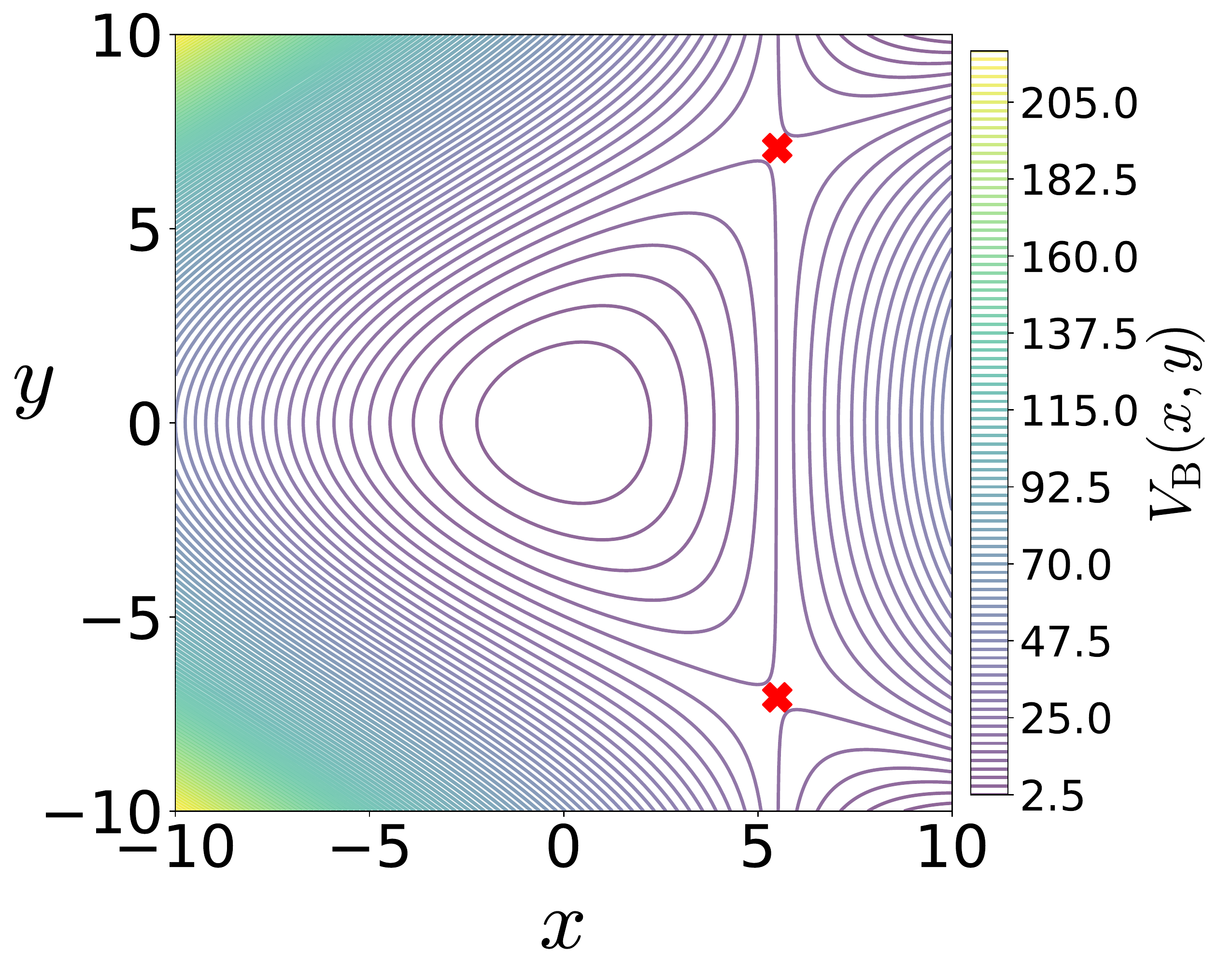}\label{fig:pes_cont_Barbanis}
	\includegraphics[width=0.25\textwidth]{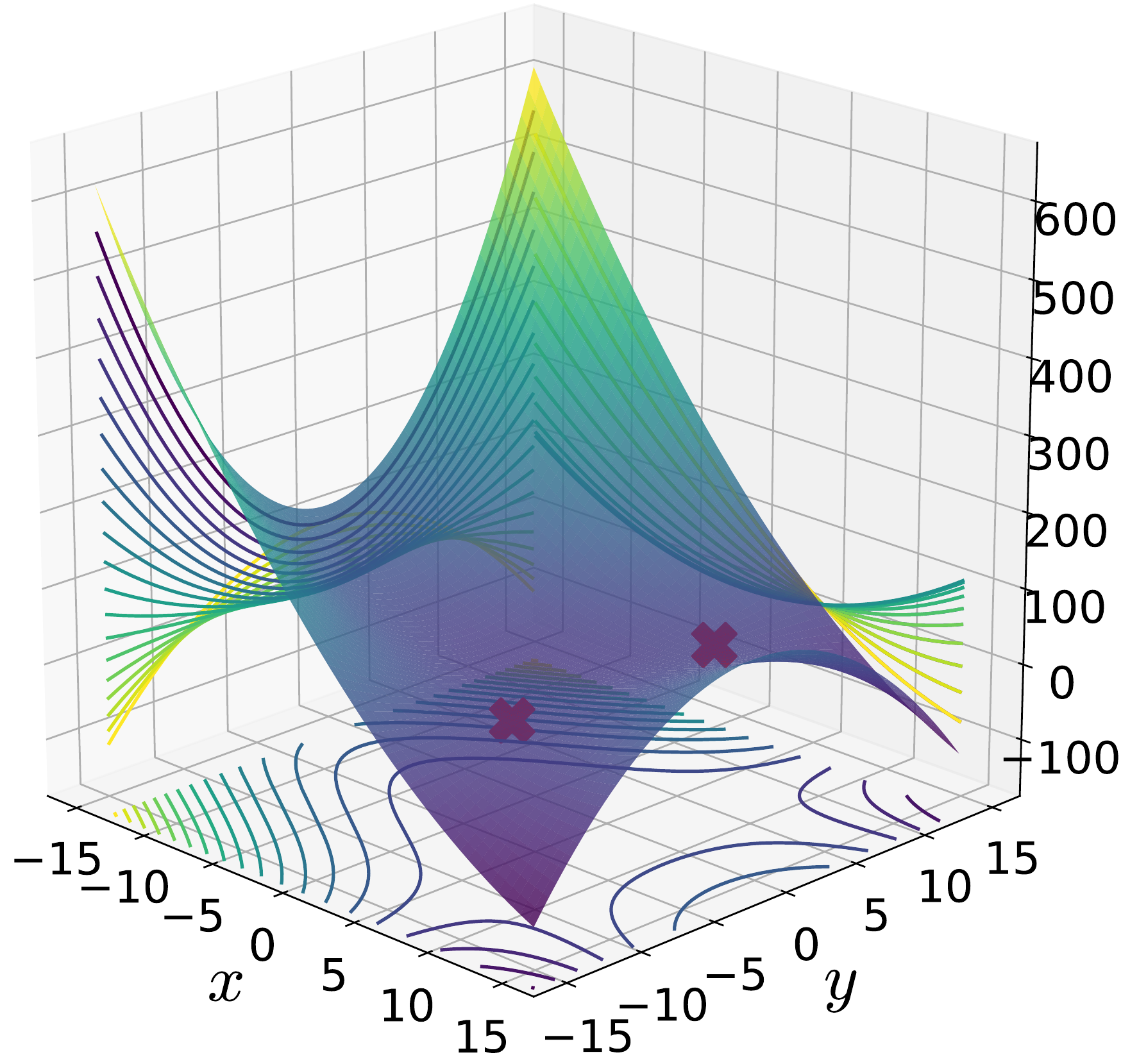}\label{fig:pes_cont_Barbanis_3D}}
	\subfigure[]{\includegraphics[width=0.4\textwidth]{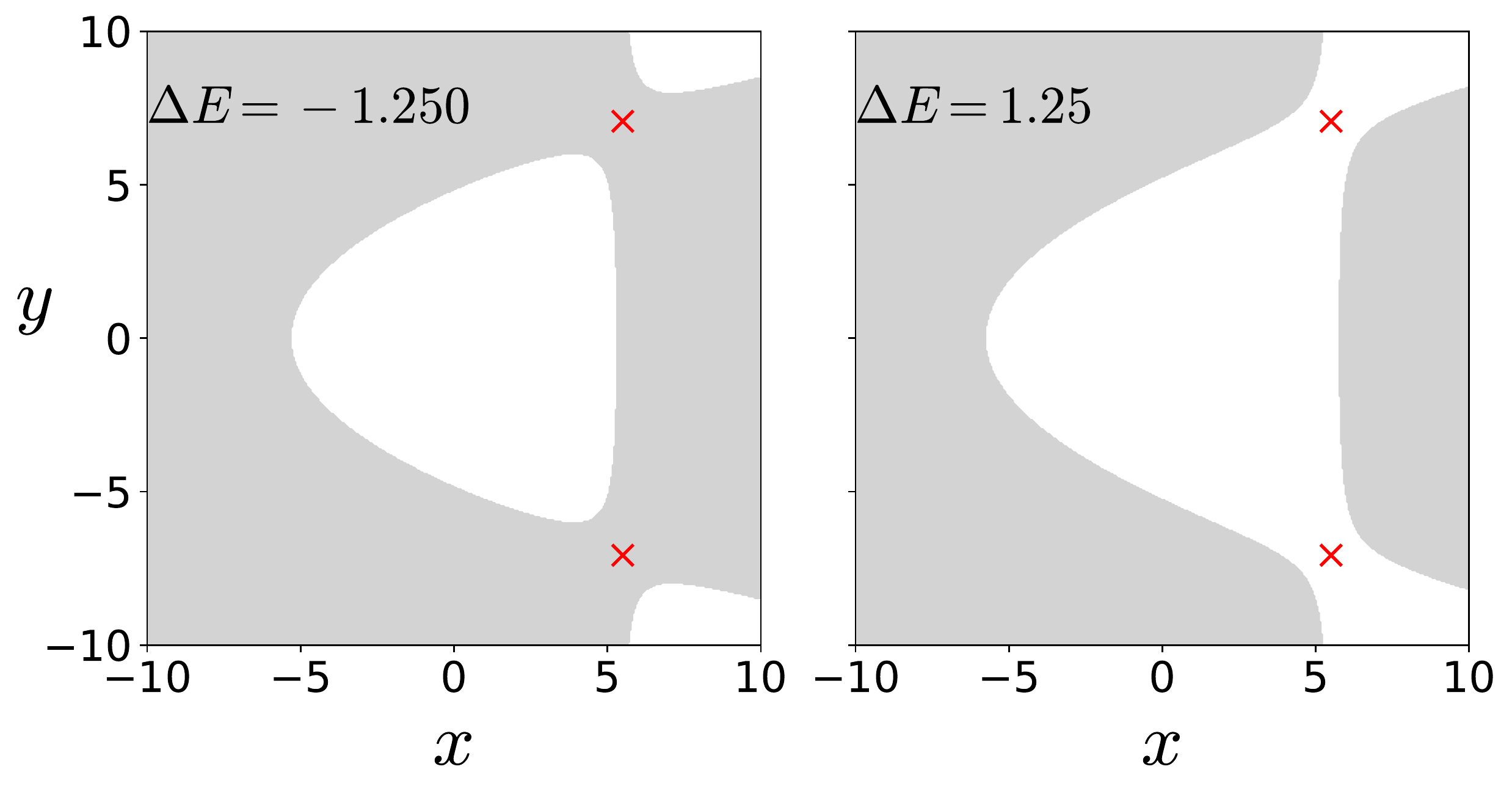}\label{fig:hills_region_Barbanis}}
	\caption{(a) Potential energy function underlying the coupled harmonic Hamiltonian~\eqref{eqn:Hamiltonian_Barbanis} as isopotential contour and surface. The index-1 saddles are shown as red crosses in both the plots. (b) Hill's region for energy below and above the energy of the index-1 saddle. Parameters used are 
	$\omega_x = 1.0, \omega_y = 1.1, \delta = -0.11$.}
\end{figure*}

The contours of the coupled harmonic 2 DoF potential energy function in 
~\eqref{eqn:Hamiltonian_Barbanis} is shown in Fig.~\ref{fig:pes_cont_Barbanis} along with the 3D 
view of the surface. We note here that 
the potential has steep walls for $x < 0$ when $\delta < 0$ and steep drop-off beyond the 
bottlenecks around the index-1 saddles. This leads to unphysical motion in the sense of 
trajectories approaching $-\infty$ with ever increasing acceleration even for finite values of the 
configuration space coordinates~\cite{brumer_variational_1976}. 

In Fig.~\ref{fig:hills_region_Barbanis} we show the \emph{Hill's region}, as defined in 
App.~\ref{ssect:coupled_2dof}, for the model system~\eqref{eqn:Hamiltonian_Barbanis}. It is important to note here that even 
though Hill's region is shown on the configuration space, it captures the dynamical 
picture, that is the {\em phase space perspective}, of the Hamiltonian. This visualization of the 
energetically accessible and forbidden realm is the first step towards introducing 
two-dimensional surfaces to explore trajectory behavior. The complete description of the unstable periodic orbit and its invariant manifolds is described in 
App.~\ref{ssect:tube_mani} along with the visualization in the 3D space.

% %height=2in
%

% \textbf{Low dimensional probe of the phase space structures~\textemdash~}
%\textit{Intersection with isoenergetic two-dimensional surface~\textemdash~} 
Since this model system is conservative 2 DoF Hamiltonian, that is the phase space is $\mathbb{R}^4$, the energy surface is three dimensional, 
the dividing surface is two dimensional, and the normally hyperbolic invariant manifold (NHIM), referred to as the unstable periodic orbit, is one dimensional~\cite{wiggins_role_2016}. 
%With this in mind, let $\mathcal{M}(e)$ be the energy surface given by setting the 
%Hamiltonian~\eqref{eqn:Hamiltonian_Barbanis} equal to a 
%constant $e$, that is
%% 
%\begin{align}
%	\mathcal{M}(e) = \{(x,y,p_x,p_y) \; | \; \mathcal{H}(x,y,p_x,p_y) = e \} \label{eqn:energy_surf}
%\end{align}
%
%and is a three-dimensional surface embedded to obtain a low dimensional sample of the dynamics.
Now, if we consider the intersection of a two dimensional surface 
with the three-dimensional energy surface, we would obtain the one-dimensional energy boundary on 
the surface of section\footnotetext{Dimension of intersection of object 1 and object 2 = 
Dimension of object 1 + Dimension of 
object 2 - Dimension of ambient space. This dimensional argument holds for all surfaces as along as 
they are not tangential or coincide with each other}. We will focus our study by using the 
isoenergetic two-dimensional surface
\begin{align}
%	 U_{V, +} &= \left\{(y,p_y) \; | \; x = 0, \; p_x(y,p_y;e) > 0 \right\} , \qquad \text{motion 
%to the right}  \label{eqn:sos_U_Vp} \\
%    U_{V, -} &= \left\{(y,p_y) \; | \; x = 0, \; p_x(y,p_y;e) < 0 \right\} , \qquad \text{motion 
%to the left}  \label{eqn:sos_U_Vm} \\
	U_{xp_x, +} &= \left\{(x,y,p_x,p_y) \; | \; y = 0, \; p_y(x,y,p_x;e) > 0 \right\} 
	\label{eqn:sos_Uxpx} 
%    U_{H, -} &= \left\{(x,p_x) \; | \; y = 0, \; p_y(x,p_x;e) < 0 \right\} , \qquad \text{motion 
%to the bottom}  \label{eqn:sos_U_Hm}
\end{align}
where the sign of the momentum coordinate enforces a directional crossing of the surface. Due to 
the form of the vector field~\eqref{eqn:two_dof_Barbanis} and choice of $\delta < 0$, this 
directionality condition implies motion towards positive $y$-coordinate. 

In this article, detecting the phase space structures will constitute finding the intersection of the NHIM and its invariant manifolds with a two dimensional surface (for example, Eqn.~\eqref{eqn:sos_Uxpx}).

% This will be achieved using a trajectory diagnostic on this surface.
%For finding the NHIM and its invariant manifolds, we will use the 
%surface-of-section~\eqref{eqn:sos_Uxpx} for constructing Poincar\'e return map and 
%elucidating the advantages of Lagrangian descriptor. 
%One approach to compute this SOS is to initialize points on the energy surface, $E(x,y,p_x,p_y) = 
%e$, along the 
%x-coordinate.  

\subsection{Model system: coupled harmonic 3 DoF Hamiltonian}
\label{sec:model_prob_3dof}

%has been discussed in the literature under the name of for 3 degrees-of-freedom galactic model
The next higher dimensional model system to consider is the coupled harmonic potential in 3 
dimensions and underlying a 3 degrees-of-freedom system in~\cite{contopoulos_1994,farantos_1998}. 
The Hamiltonian is given by
%
% \begin{widetext}
\begin{equation}
\begin{split}
\mathcal{H}(x,y,z,p_x,p_y,p_z) =  T(p_x, p_y, p_z) + V_{\rm BC}(x,y,z) \\ 
\mathrel{\phantom{=}} =  \frac{1}{2}p_x^2 + \frac{1}{2}p_y^2 + \frac{1}{2}p_z^2  + \\ \frac{1}{2}\omega_x^2 x^2 + \frac{1}{2}\omega_y^2 y^2 + \frac{1}{2}\omega_z^2 z^2 - \epsilon x^2y - \eta x^2 z  
\label{eqn:Hamiltonian_BC_3dof}
\end{split}
\end{equation}
% \end{widetext}

where $\omega_x^2, \omega_y^2, \omega_z^2, \epsilon, \eta$ are the parameters related to the 
coupled harmonic 3 dimensional potential energy function~\cite{farantos_1998}. In this study, we 
will fix the parameters to be $\omega_x^2 = 0.9, \omega_y^2 = 1.6, \omega_z^2 = 0.4, \epsilon = 
0.08, \eta = 0.01$. The two index-1 saddle equilibria (as shown in the 
App.~\ref{sect:coupled_3dof}) of the Hamiltonian vector field~\eqref{eqn:three_dof_Barbanis} 
are located at 
%
% \begin{widetext} \mathbf{x}_{\rm eq} = 
\begin{equation}
\begin{split}
\left(\pm \frac{\omega_x\omega_y\omega_z}{\sqrt{2(\epsilon^2\omega_z^2 + 
		\eta^2\omega_y^2)}}, 
\frac{\epsilon \omega_x^2\omega_z^2}{2(\epsilon^2\omega_z^2 + \eta^2\omega_y^2)}, \right. \nonumber \\
\left.\kern - \nulldelimiterspace \frac{\eta \omega_x^2\omega_y^2}{2(\epsilon^2\omega_z^2 + \eta^2\omega_y^2)},  0, 0, 0 \right) 
\label{eqn:eq_pt_BC_3dof}
\end{split}
\end{equation}
% \end{widetext}
and the total energy is 
\begin{equation}
E_c = \frac{1}{8} \omega_x^2 \frac{\omega_x^2 \omega_y^2 \omega_z^2}{ \left( \epsilon^2 \omega_z^2 
	+ \eta^2 \omega_y^2 \right)}.
\end{equation}
The equilibrium point at $(0,0,0,0,0,0)$ is stable and has total energy $0$. For the parameters 
used in this study, the equilibrium points are located at $\left( \pm 10.290, 5.294, 2.647, 
0, 0, 0 \right)$ and $\left( 0, 0, 0, 0, 0, 0 \right)$ and have total energy, $E_c \approx 
23.824$ and $E = 0$, respectively.
% 
%\begin{equation}
%\left( 10.28992, 5.29412, 2.64706, 0, 0, 0 \right) \qquad \text{and} \qquad \left( 0, 0, 0, 0, 0, 
%0 \right)
%\end{equation}
%

\begin{figure*}[!th]
	\centering
	\subfigure[]{\includegraphics[width=0.3\textwidth]{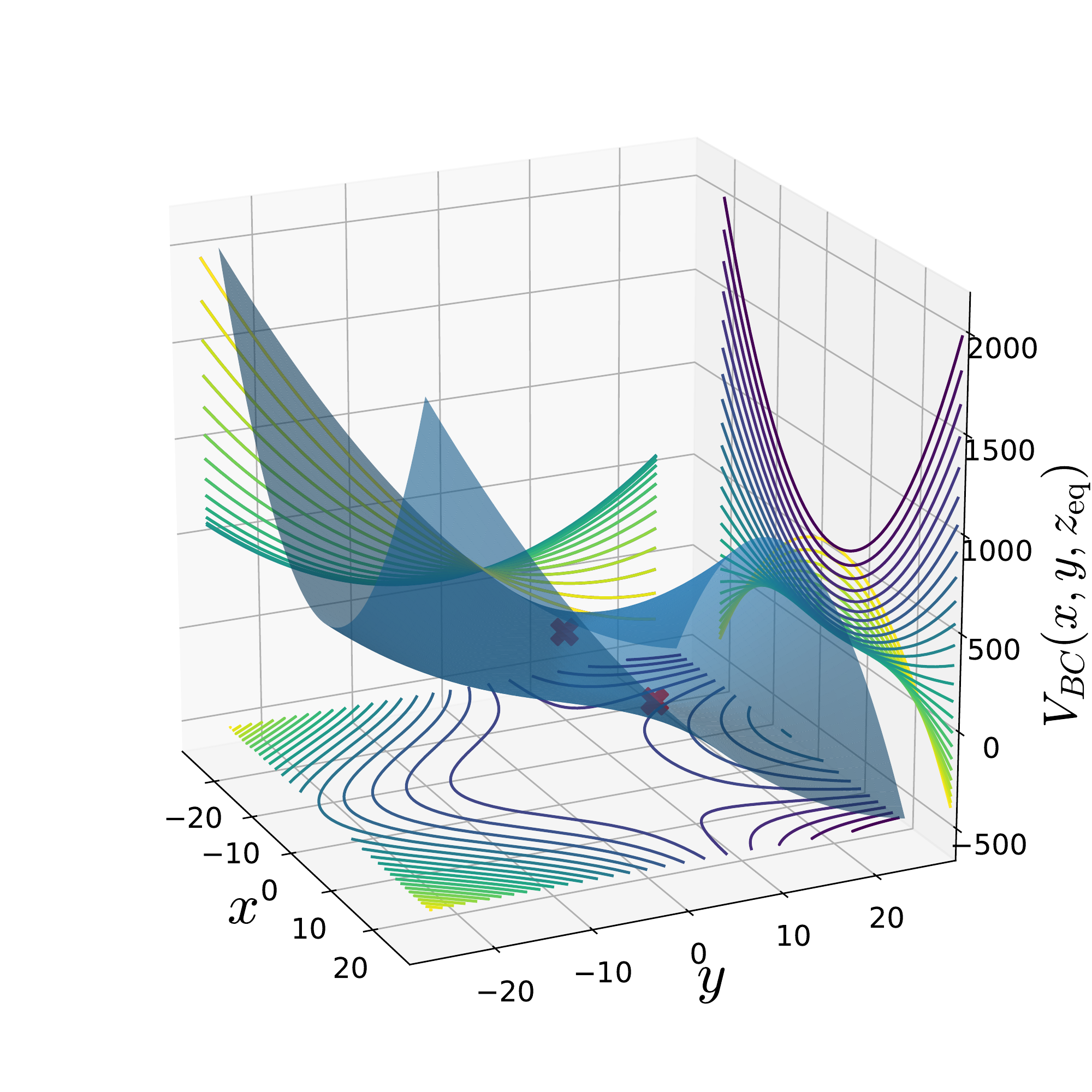}} 
	%width=0.35\textwidth
	\subfigure[]{\includegraphics[height=1.5in]{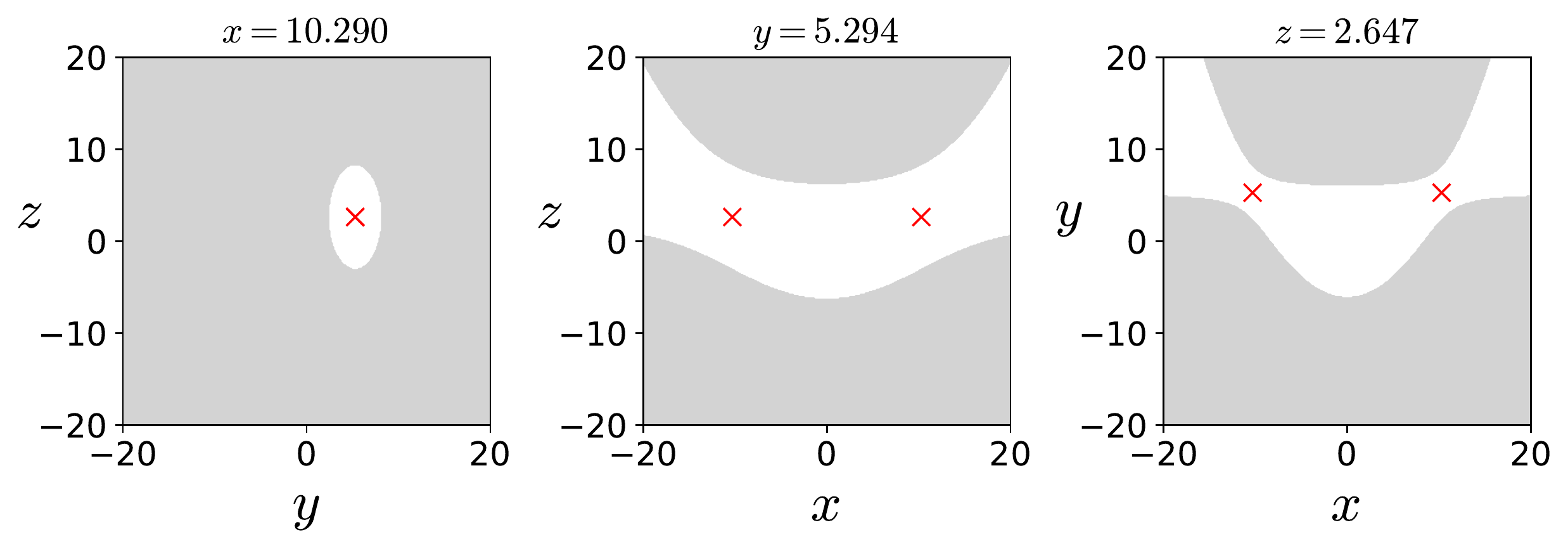}} 
	%width=0.6\textwidth
%	\includegraphics[width=0.85\textwidth]{./figures/Barbanis3dof_pes2d300x300_zslices.png}
%	\caption{\textbf{Contour plot} of the Barbanis-Contopoulos potential at $z = -50, 0, 50$ and 
%the red cross mark the critical points given by Eqn.~\eqref{eqn:eq_pt_BC_3dof}. Parameters used 
%are 
%$\omega_x^2 = 0.9, \omega_y^2 = 1.6, \omega_z^2 = 0.4, \epsilon = 0.08, \eta = 0.01$.}
	\caption{(a) Potential energy function underlying the coupled harmonic Hamiltonian~\eqref{eqn:Hamiltonian_BC_3dof} at $z_{\rm eq} = 2.647$ as isopotential contour and surface. (b) Hill's 
	region for excess energy, $\Delta E = 6.000$ and projected on the configuration space 
	coordinates at the equilibrium point. We note here that the potential energy surface and the 
	Hill's region is plotted by fixing one of the configuration coordinates at the equilibrium point.}
	\label{fig:Barbanis_Contopoulos_3dof}  
\end{figure*}

We show the isopotential contours of the potential energy function at fixed value of $z_{\rm eq}$ 
in Fig.~\ref{fig:Barbanis_Contopoulos_3dof} along with the Hill's regions for positive excess 
energy, $\Delta E = 6.000$ and projected on the configuration space coordinates at 
the equilibrium point. 

Since this model system is conservative 3 DoF Hamiltonian, that is the phase space is $\mathbb{R}^6$, the energy surface is five dimensional, the dividing surface is four dimensional, and the normally hyperbolic invariant manifold (NHIM) is three dimensional, or precisely 3-sphere, and its invariant manifolds are four dimensional, or precisely $\mathbb{R}^1 \times \mathbb{S}^3$ or {\it spherical cylinders}~\cite{wiggins_role_2016}. Now, if we consider the intersection of a two-dimensional section with the five dimensional energy surface in $\mathbb{R}^6$, we would obtain the one-dimensional energy boundary on 
the surface. We will focus our study near the bottleneck by considering the isoenergetic two dimensional surfaces 
\begin{widetext}
	\begin{align}
	U_{xp_x}^+ = & \left\{ (x, y, z, p_x, p_y, p_z) \; | \; y = y_{\rm eq}, z = z_{\rm eq}, \; p_y = 0, \; p_z(x, y, z, p_x, p_y; e) > 0 \right\} \label{eqn:Barbanis3dof_uxpx}\\
	U_{yp_y}^+ = & \left\{ (x, y, z, p_x, p_y, p_z) \; | \; x = x_{\rm eq}, z = z_{\rm eq}, \; p_x = 0, \; p_z(x, y, z, p_x, p_y; e) > 0 \right\} \label{eqn:Barbanis3dof_uypy}\\
	U_{zp_z}^+ = & \left\{ (x, y, z, p_x, p_y, p_z) \; | \; x = x_{\rm eq}, y = y_{\rm eq}, \; p_x = 0, \; p_y(x, y, z, p_x, p_z; e) > 0 \right\} \label{eqn:Barbanis3dof_uzpz}
%	\label{eqn:Barbanis3dof_sos_near_saddle}
	\end{align}
\end{widetext}

In this 3 DoF system, detecting points on the three dimensional NHIM and four dimensional invariant manifolds will constitute finding their intersection with the above two dimensional surfaces.

\subsection{\label{sec:LD}Method: Lagrangian descriptor}

The Lagrangian descriptor (LD) as presented in Ref.\cite{madrid2009} is the arc length of a 
trajectory calculated on a chosen initial time $t_0$ and measured for fixed forward and backward 
integration time, $\tau$. For continuous time dynamical systems, Ref.\cite{lopesino2017} gives an 
alternative definition of the LD which is useful for proving rigorous results and can be computed 
along with the trajectory. It provides a characterization of the notion of singular features of the 
LD that facilitates a proof for detecting invariant manifolds in certain model situations.  In 
addition, the ``additive nature'' of this new definition of LD provides an approach for assessing 
the influence of each degree-of-freedom separately.  This property was used in 
Ref.\cite{demian2017} to show that Lagrangian descriptor can detect Lyapunov periodic orbits in the 
two degrees-of-freedom H{\'e}non-Heiles system. We will adopt a similar strategy for the aforementioned two and three degrees-of-freedom autonomous H{\'e}non-Heiles type systems. 

In the general setting of a 
time-dependent vector field 
\begin{equation}
\frac{d\mathbf{x}}{dt} = \mathbf{v}(\mathbf{x},t), \quad \mathbf{x} \in \mathbb{R}^n \;,\; t \in 
\mathbb{R}
\end{equation}
where $\mathbf{v}(\mathbf{x},t) \in C^r$ ($r \geq 1$) in $\mathbf{x}$ and continuous in time. The 
definition of LDs depends on the initial condition $\mathbf{x}_{0} = \mathbf{x}(t_0)$, on the 
initial time $t_0$ (trivial for autonomous systems) and the integration time $\tau$, and the type 
of norm of the trajectory's components, and takes the form,

\begin{equation}
M_p(\mathbf{x}_{0},t_0,\tau) = \displaystyle{\int^{t_0+\tau}_{t_0-\tau} \sum_{i=1}^{n} 
|\dot{x}_{i}(t;\mathbf{x}_{0})|^p \; dt}
\label{eqn:M_function}
\end{equation}

 where $p \in (0,1]$ and $\tau \in \mathbb{R}^{+}$ are freely chosen parameters,  and the 
overdot symbol represents the derivative with respect to time. It is to be noted here that there 
are three formulations of the function $M_p$ in the literature: the arc length of a trajectory in 
phase space~\cite{madrid2009}, the arc length of a trajectory projected on the configuration 
space~~\cite{junginger2016transition,junginger2016uncovering,junginger2017chemical,junginger2017variational},
 and the sum of the $p$-norm of the vector field components~\cite{lopesino_2015,lopesino2017}.
% In the context of chemical reactions, a different form with $p = 1$ has been popular in 
%Refs.~\cite{junginger2016transition,junginger2016uncovering,junginger2017chemical,junginger2017variational}.(I
% thought these people used the arclength formulation of LD and not the definition of Lopesino?) 
Although the latter formulation of the Lagrangian descriptor~\eqref{eqn:M_function} developed in 
Ref.~\cite{lopesino_2015,lopesino2017} does not resemble the arc length, the numerical results 
using either of these forms have been shown to be in agreement and promise of predictive capability 
in geophysical flows~\cite{amism11,mendoza2014,ggmwm15,ramos2018}. The formulation we adopt here is 
motivated by the fact that this allows for proving rigorous result, which we will discuss in the 
next section, connecting the singular features and minimum in the LD plots with NHIM and its stable 
and unstable manifolds. 
% In a forthcoming paper on{PROOF (what does this mean?): singular and extremal values in LD detect 
%manifolds}, we present the rigorous result from which it follows that {\color{blue} SN: referring 
%to the manuscript that we are preparing with Victor.}
It follows from the result that 
\begin{align}
\mathcal{W}^s(\mathbf{x}_0, t_0) & = \text{\rm argmin} \; \mathcal{L}^{(f)}(\mathbf{x}_0, t_0, 
\tau) \\
\mathcal{W}^u(\mathbf{x}_0, t_0) & = \text{\rm argmin} \; \mathcal{L}^{(b)}(\mathbf{x}_0, t_0, \tau)
\end{align}
where the stable and unstable manifolds ($\mathcal{W}^s(\mathbf{x}_0, t_0)$ and 
$\mathcal{W}^u(\mathbf{x}_0, t_0)$) denote the invariant manifolds at intial time $t_0$ and 
$\text{\rm argmin} \; (\cdot)$ denotes the argument that minimizes the function 
$\mathcal{L}^{(\cdot)}(\mathbf{x}_0, t_0, \tau)$ in forward and backward time, respectively. In 
addition, the coordinates on the NHIM, $\mathcal{M}(\mathbf{x}_0, t_0)$ at time $t_0$ is given by the intersection 
$\mathcal{W}^s(\mathbf{x}_0, t_0)$ and $\mathcal{W}^u(\mathbf{x}_0, t_0)$ of the stable and 
unstable manifolds, and thus given by
\begin{equation}
\begin{aligned}
\mathcal{M}(\mathbf{x}_0, t_0) & = \text{\rm argmin} \; \left( \mathcal{L}^{(f)}(\mathbf{x}_0, t_0, 
\tau) + \mathcal{L}^{(b)}(\mathbf{x}_0, t_0, \tau) \right) \\
& = \text{\rm argmin} \; 
\mathcal{L}(\mathbf{x}_0, t_0, \tau) 
\end{aligned}
\end{equation}

% \textbf{Why we need to modify fixed integration time LD~\textemdash~}

In applying the LD method to nonlinear systems, one observes multiple minima and singularities that can lead to trouble with isolating the one minima due to the NHIM and the ones due to its invariant manifolds. 
Since, as we integrate initial conditions on an isoenergetic two dimensional surface such as  $U_{xp_x}^+$~\eqref{eqn:sos_Uxpx}, almost all 
trajectories that escape to infinity get integrated for the entire time interval and result in 
numerical overflow of the function M value~\eqref{eqn:M_function} and show up as NaN. This can, however, be avoided by integrating for shorter time interval but this will vary for different locations of a surface. Thus, leading to trouble in locating the point with minimum and singularity in LD contour map that correspond to NHIM and its invariant manifolds.

\begin{figure*}[!ht]
	\centering
	\subfigure{\includegraphics[width=0.85\textwidth]{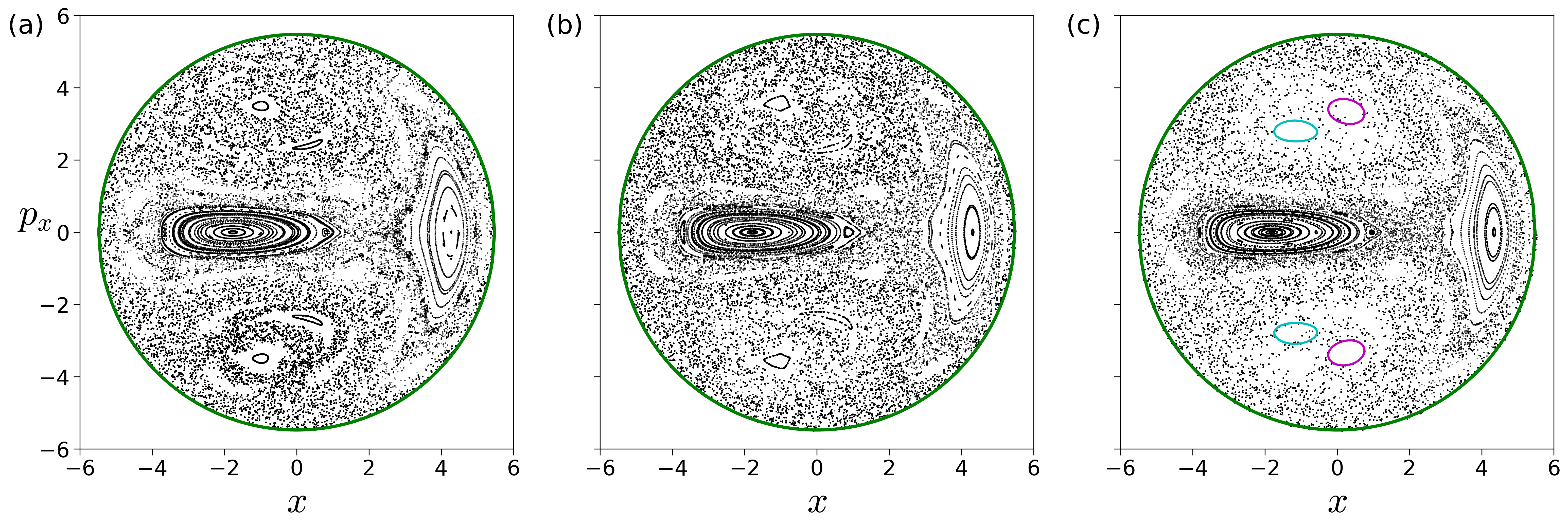}}
    \subfigure{\includegraphics[width=0.85\textwidth]{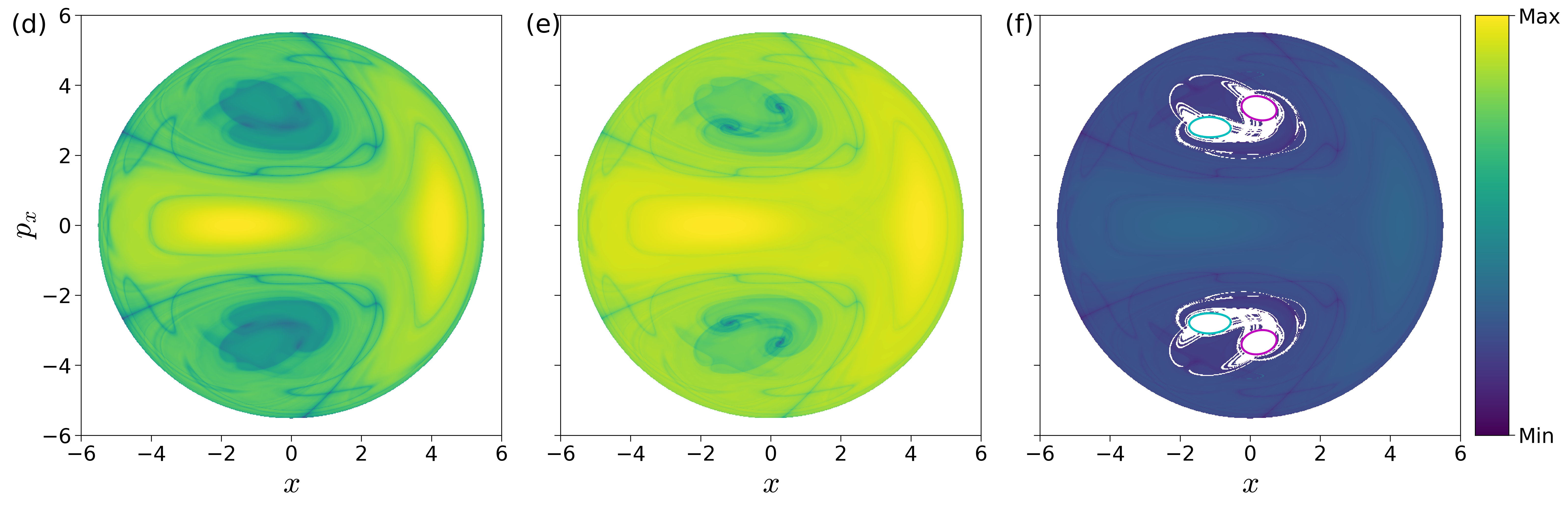}}
% 	\subfigure[]{\includegraphics[width=2in]{./figures/lag_desc_Barbanis_2dof_mani_sos_E15-250_tau50.png}\label{fig:lag_desc_Barbanis_E15-250_x-px_tau50_scaled-arctan}}
    \caption{\textbf{Top row:} Poincar\'e surface of section, $U_{xp_x}$~\eqref{eqn:sos_Uxpx}, at excess energy (a) 
	$\Delta E = -0.125$, (b) $\Delta E = 0.000$, (c) $\Delta E = 0.125$ where the intersection of the surface of section with the energy surface is shown in green. \textbf{Bottom row:} Lagrangian descriptor on the surface of section, $U_{xp_x}$~\eqref{eqn:sos_Uxpx}, for the excess energies (d) 
	$\Delta E = -0.125$, (e) $\Delta E = 0.000$, (f) $\Delta E = 0.125$ and the integration time 
	$\tau = 50$. The intersection of the surface of section with the cylindrical manifolds of the NHIM \textemdash unstable periodic orbit for this system \textemdash associated with the index-1 saddle equilibrium point in the bottleneck is shown in cyan (stable) and magenta (unstable) curves. The magenta and cyan curves in $p_x > 0$ correspond to the invariant manifolds of unstable periodic orbit at $y > 0$ index-1 saddle, and the ones in $p_x < 0$ correspnd to the invariant manifolds of unstable periodic orbit at $y < 0$ index-1 saddle.} 
	\label{fig:psect_lag_desc_Barbanis}
\end{figure*}

This computational issue has been addressed in recent efforts to locate transition state trajectory in driven and 3 degrees-of-freedom chemical reaction dynamics~\cite{craven2017lagrangian,craven2016deconstructing,craven2015lagrangian}. It has been noted that computing fixed integration time Lagrangian descriptor (LD) leads to two 
potential issues:

1. Bounded trajectories  will show global recrossings of the barrier as predicted by Poincar\'e 
recurrrence theorem. The recrossings will show multiple minima and singularities (as in Fig.~\ref{fig:psect_lag_desc_Barbanis}(d-f)) in the LD plot 
which obscures locating the actual NHIM. 

2. The trajectories that escape the potential well will leave with ever increasing acceleration, if 
the potential energy surface opens out to infinity. The trajectories with NaN LD values will 
render the contour map flat which again obscures locating the NHIM. 

To circumvent these issues, a heuristic that has been adopted in the literature is to calculate LD 
values only until a trajectory remains inside the barrier region. The immediate result is the initial condition on an invariant manifold will have a maxima in the LD values because of being integrated for the full integration time (preselected) interval. 

Thus, the formulation~\eqref{eqn:M_function} is modified as
\begin{equation}
M_p(\mathbf{x}_{0},t_0,\tau^{\pm}) = \displaystyle{\int^{t_0+\tau^+}_{t_0-\tau^-} \sum_{i=1}^{n} 
	|\dot{x}_{i}(t;\mathbf{x}_{0})|^p \; dt}
\label{eqn:M_function_var}
\end{equation}
where the integration time interval depends on a trajectory and given by
\begin{equation}
	\tau^{\pm}(\mathbf{x}_0) = \min\left(\tau, t\vert_{|\mathbf{x}(t)| > q_s} \right)
	\label{eqn:var_time_qs}
\end{equation}
where $q_s$ defines a domain, called the {\it saddle region}, in the configuration space around the saddle. We note here that the only initial condition that gets integrated for the entrire $\tau$ time units in forward and backward time is the one on the NHIM. In 
addition, the coordinates on the NHIM, $\mathcal{M}(\mathbf{x}_0, t_0)$, at time $t_0$ is given by
\begin{equation}
\begin{aligned}
\mathcal{M}(\mathbf{x}_0, t_0) & = \text{\rm argmax} \; \left( \mathcal{L}^{(f)}(\mathbf{x}_0, t_0, 
\tau) + \mathcal{L}^{(b)}(\mathbf{x}_0, t_0, \tau) \right) \\
& = \text{\rm argmax} \; \mathcal{L}(\mathbf{x}_0, t_0, \tau) 
\end{aligned}
\end{equation}
This is also a familiar from a dynamical systems perspective where the literature on average exit times to locate invariant sets has been discussed for the symplectic maps (see~\cite{meiss_average_1997} and related references). However, the connection between features in exit times and LD contour maps is not the focus of this study and will be deferred as related future work.

\section{\label{sec:results}Results}

% \subsection{Detecting Phase Space Structures using the Lagrangian Descriptor}

We begin by noting that two-dimensional Poincar\'e surface of section have sufficient 
dimensionality to capture trajectories on a three dimensional energy surface, however for high 
dimensional systems trajectories can go ``around'' the two dimensional surface. 
One approach available in the literature is to use high dimensional Poincar\'e sections which can 
``catch'' trajectories but are hard to visualize on paper or in the virtual 3D space. Even when 
gets around this issue, using suitable projective geometry, the fact that the qualitative analysis 
based on Poincar\'e sections depends on trajectories returning to this surface can not be 
circumvented since trajectories on and inside the spherical cylinders will not return to the 
Poincar\'e surface of section.

\subsection{Coupled harmonic 2 DoF system} 

As discussed in aforementioned literature~\cite{madrid2009,lopesino2017,demian2017}, points with 
minimum Lagrangian descriptor (LD) values and singularity are on the invariant manifolds. In 
addition, LD plots show dynamical correspondence with Poincar\'e sections (in the sense that 
regions with regular and chaotic dynamics are distinct in both Poincar\'e section and LD plots) 
while also depicting the geometry of manifold 
intersections~\cite{demian2017,lopesino2017,garcia-garrido_2018}. This correspondence in the LD 
features and Poincar\'e section is confirmed in Fig.~\ref{fig:psect_lag_desc_Barbanis} where we 
show the Poincar\'e surface of section Eqn.~\eqref{eqn:sos_Uxpx} of trajectories and LD contour 
maps on 
the same isoenergetic two-dimensional surface for negative and positive excess energies. It can be 
seen that the chaotic dynamics as marked by the sea of points in Poincar\'e section is revealed as 
the tangle of invariant manifolds which are points of minima and singularity in the LD plots. As 
shown by the one dimensional slices of the LD plots, there are multiple such minima and 
singularities and as the excess energy is increased to positive 
values, there are regions of discontinuities along the one dimensional slice. Next, as the energy is increased and the bottleneck opens at critical energy 
$E_c$, trajectories that leave the potential well and do not return to the surface of section are 
not observed on the Poincar\'e section while the LD contour maps clearly identifies these regions 
as discontinuities in the LD values. These regions lead to escape because they are inside the 
cylindrical manifolds of the unstable periodic orbit associated with the index-1 saddle equilibrium 
point~\cite{NaRo2017}. These regions on the isoenergetic two-dimensional surface are also referred 
to as {\it reactive islands} in chemical reaction dynamics~\cite{DeLeon1992,DeMeTo1991,DeMeTo1991a}.
The escape regions or reactive islands that appear over the integration time interval can also be 
identified by using the forward and backward LD contour maps where these regions appear as 
discontinuities. In Fig.~\ref{fig:psect_lag_desc_Barbanis}(f), we show these for $\Delta E = 
0.125$ and $\tau = 50$ along with the intersection of the cylindrical manifolds' intersections that 
are computed using differential correction and numerical continuation. The detailed comparison and 
extension to high dimensional systems is not the focus of this study and will be discussed in 
forthcoming work. Thus LD maps also provide a quick and reliable approach for detecting regions 
that will lead to escape within the observed time, or in the computational context, the integration 
time. 
%\improvement[inline]{need to update the figures with all the reactive islands laid out for the 
%integration time used to obtain LD plot.}

%The color scale is normalized to the range $[500, 600]$ to elucidate the features. 

\begin{figure*}[!ht]
	\centering
	\subfigure[$U_{xp_x}^+(-7.0)$]{\includegraphics[width=0.2\textwidth]{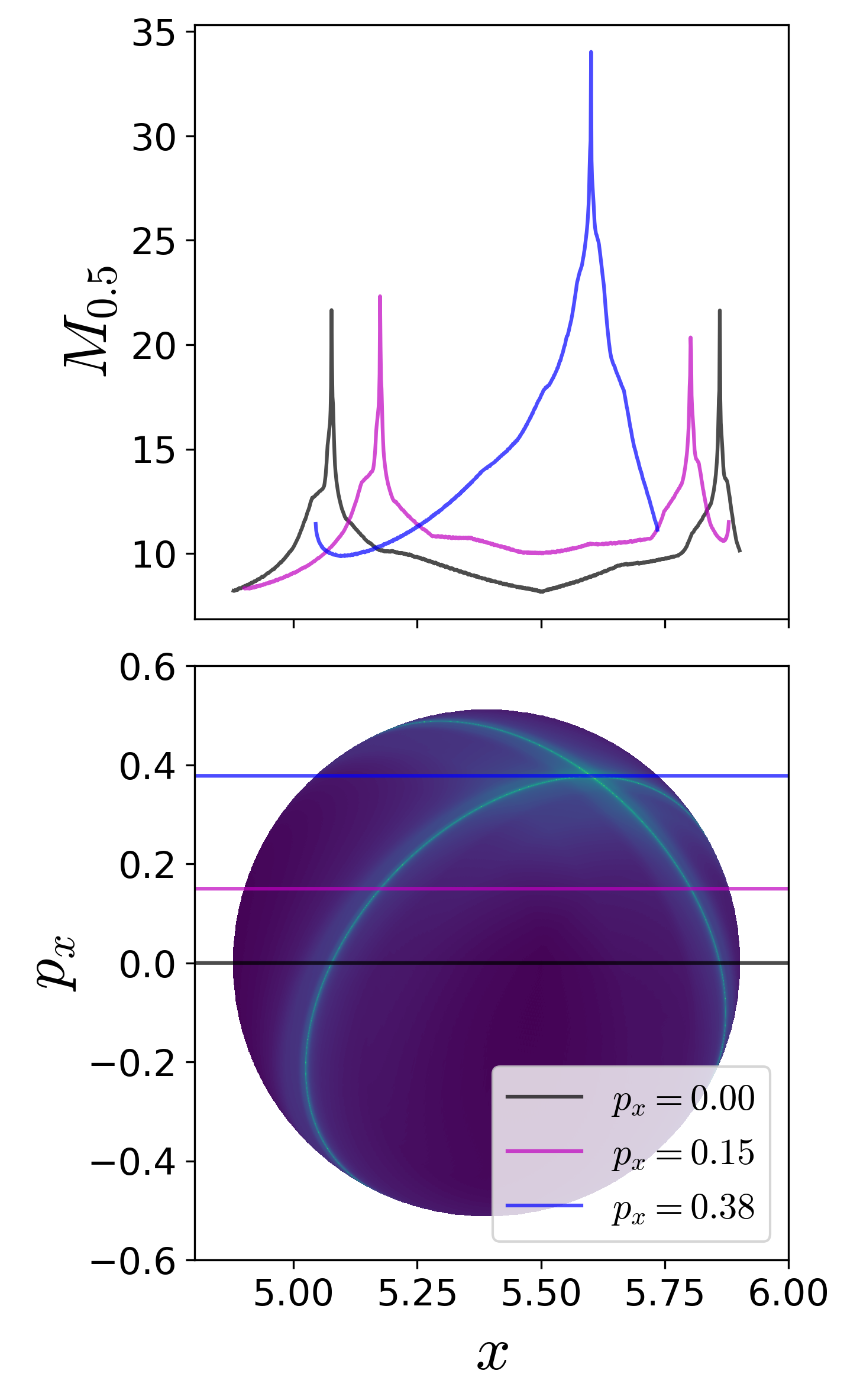}}
	\subfigure[$U_{xp_x}^+(-7.1)$]{\includegraphics[width=0.2\textwidth]{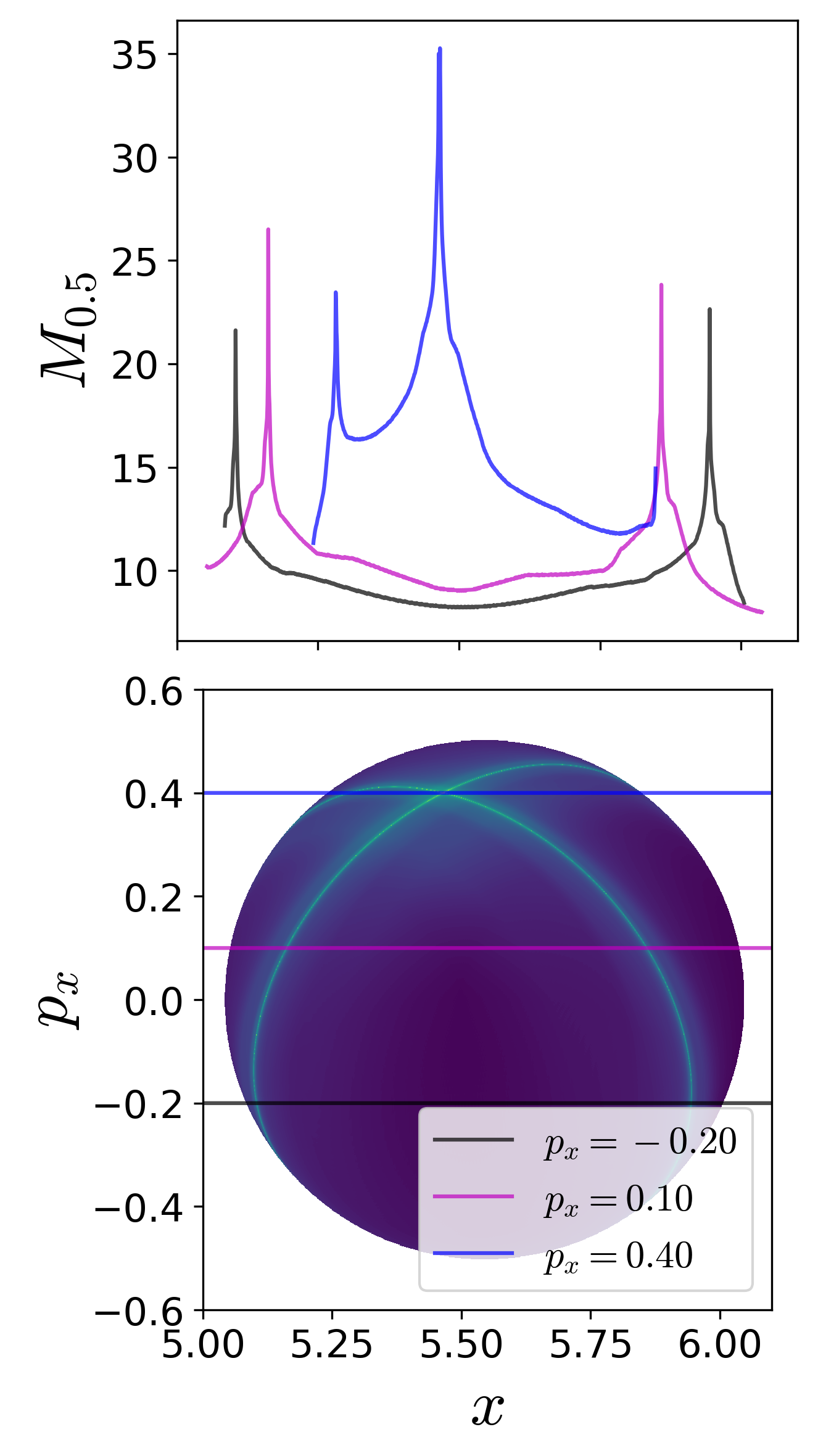}}
	\subfigure[$U_{xp_x}^+(-7.2)$]{\includegraphics[width=0.2\textwidth]{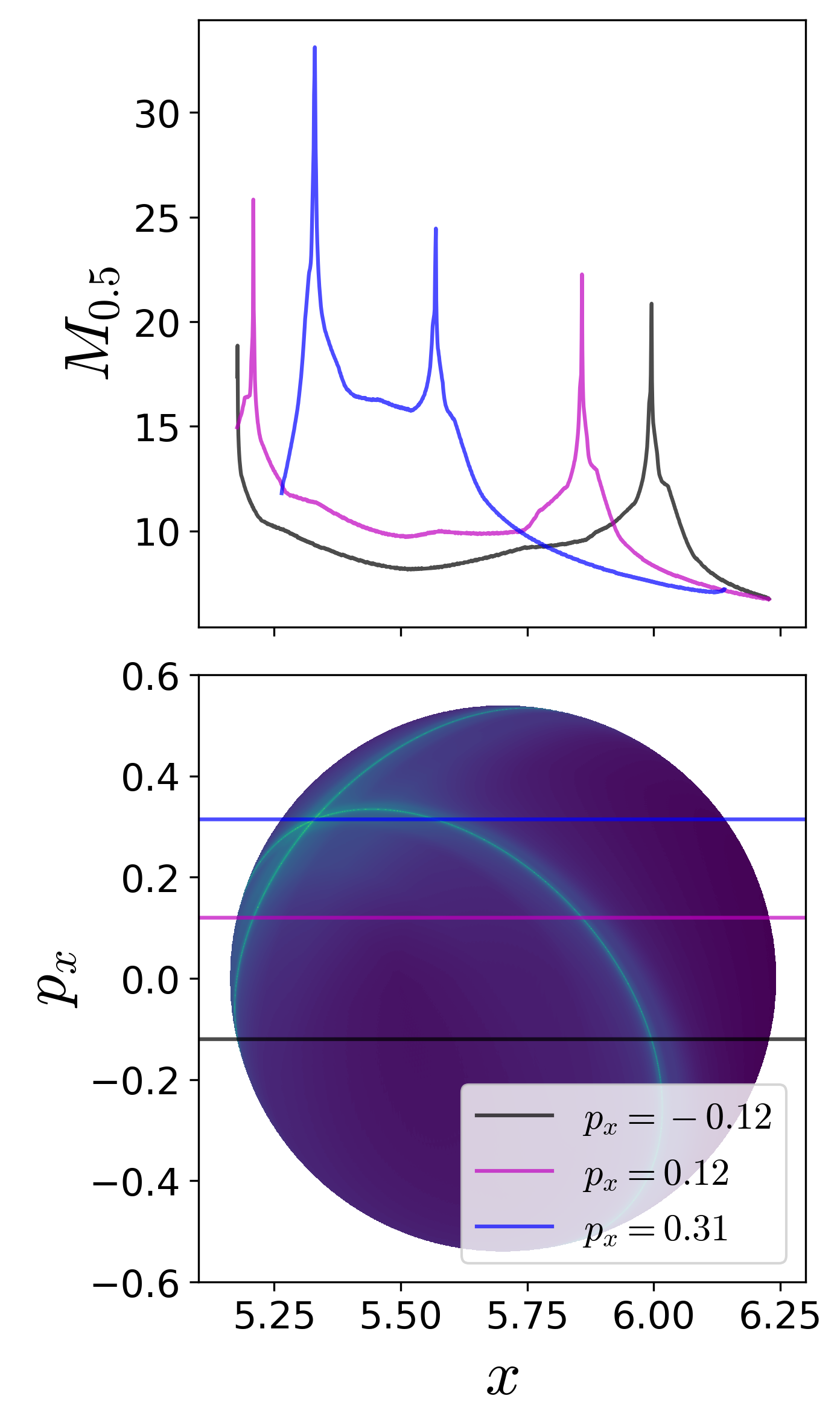}}
	\subfigure[]{\includegraphics[width=0.33\textwidth]{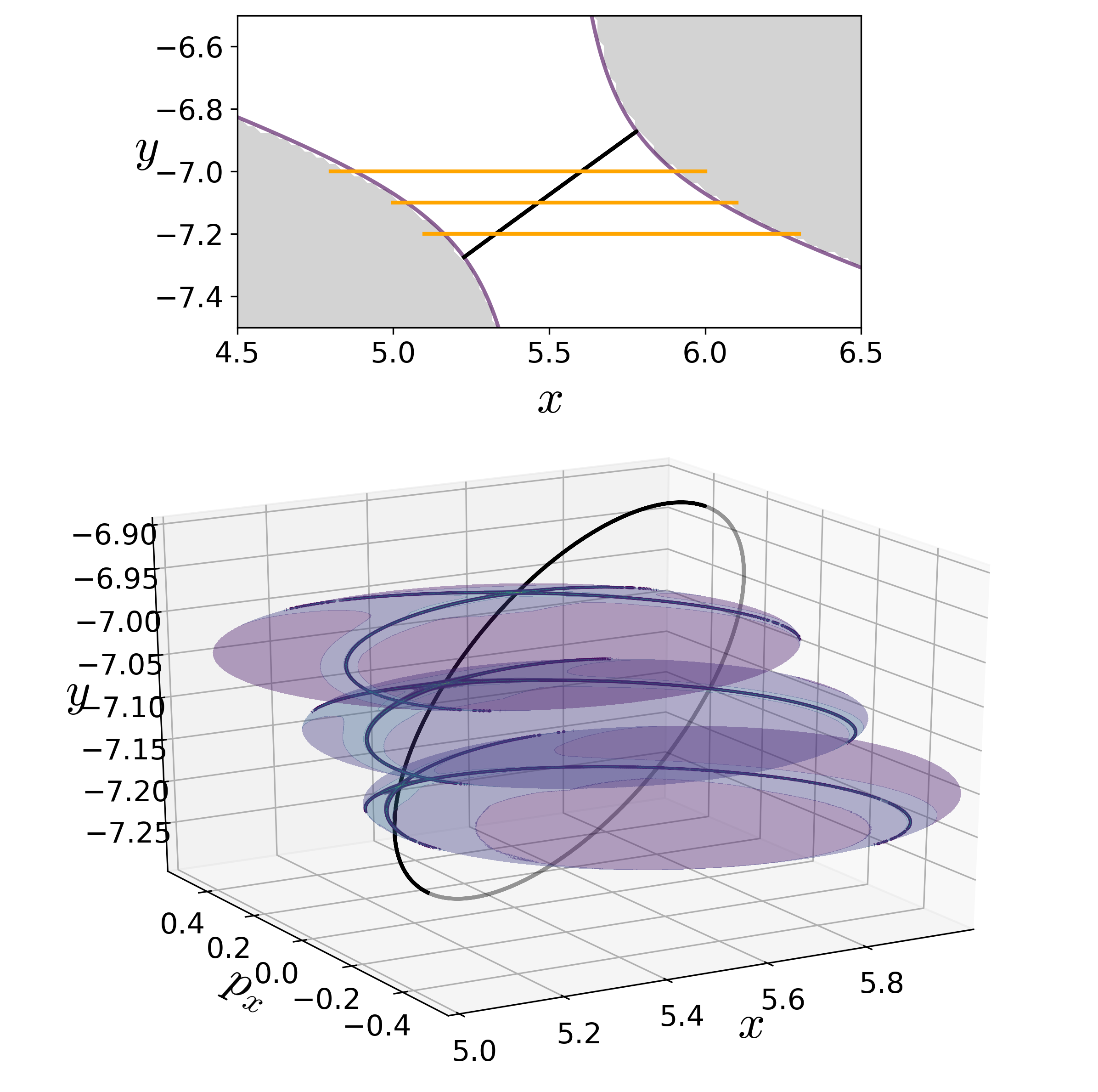}}
	\caption{Lagrangian descriptor computed for variable integration time on two dimensional slices~\eqref{eqn:sos_xpx_k} near the bottleneck that detect the NHIM and its invariant manifolds associated with the index-1 saddle. The two dimensional surfaces are shown in the top figure in (d) projected as orange lines on the configuration space and the unstable periodic orbit as black line connecting the isopotential contour corresponding to $\Delta E = 0.125$ with the Hill's region shown in grey. The two dimensional slices represent low dimensional probe of the unstable periodic orbit and the movie of a rotating view can be found \href{https://youtu.be/Ai65ULLlljc}{here}.}
	\label{fig:Barbanis2dof_M1500x1500_E15-250}
\end{figure*}

To detect the NHIM \textemdash~in this case, unstable periodic orbit \textemdash~associated with the index-1 saddles (marked by cross in Fig.~\ref{fig:pes_cont_Barbanis}), we define an isoenergetic two dimensional surface 
that is parametrized by the $y$-coordinate and placed near the $x$-coordinate of the saddle 
equilibrium that has the negative $y$-coordinate. This can be expressed as a parametric two dimensional surface 
\begin{align}
%  U_{xy}^+ = & \left\{(x,y,p_x,p_y) \; | \; p_x = 0, \; p_y(x,y;e) > 0 \right\}, \\
U_{xp_x}^+(k) = & \left\{(x,y,p_x,p_y) \, | \, y = k, p_y(x,y,p_x;e) > 0 \right\} 
\label{eqn:sos_xpx_k}
%  U_{yp_y}^+ = & \left\{(x,y,p_x,p_y) \; | \; x = k_x, \; p_y(x,y,p_x;e) > 0 \right\},
\end{align}
for total energy, $e$, which is above the critical energy, $E_c$, $k$ is the $y$-coordinate. The variable integration time LD contour maps are shown in 
Fig.~\ref{fig:Barbanis2dof_M1500x1500_E15-250} along with the projection of the low dimensional slices~\eqref{eqn:sos_xpx_k} in the configuration space and the NHIM. The points on NHIM, which is an unstable periodic orbit for 2 DoF, on this surface is the coordinate with maximum (for variable integration time) LD value. The full visualization of the NHIM as the black ellipse, $\mathbb{S}^1$, is in Fig.~\ref{fig:Barbanis2dof_M1500x1500_E15-250}(d) and has been computed using differential correction and numerical continuation (details in App.~\ref{ssect:tube_mani}) and shows clearly that points on this unstable periodic orbit are detected by the LD contour map.
%given by the intersection of the invariant manifolds 
%and also where the singular feature coincides with the global minimum of the LD values.
 %height=0.32\textheight

\subsection{Coupled harmonic 3 DoF system}
% \textbf{Detecting NHIM associated with the index-1 saddle~\textemdash~} 

The Lagrangian descriptor based approach for detecting NHIM in 2 DoF system can now be applied to 
the 3 DoF system~\eqref{eqn:Hamiltonian_BC_3dof}. On the five dimensional energy surface, the phase 
space structures such as the NHIM and its invariant manifolds are three and four 
dimensional, respectively~\cite{wiggins_role_2016}. As noted earlier, direct visualization techniques will fall short in 4 or more DoF systems even if they are successful in 2 and 3 DoF. So, LD based approach can be used to detect points on a NHIM and its invariant manifolds using low dimensional probe which are based on trajectory diagnostic on an isoenergetic two dimensional surface.  

It is to be noted that the increase in phase space dimension, leads to a polynomial scaling in the number of coordinate pairs (that is $2N(2N-1)(N-1)$ coordinate pairs for $N$ DoF system) and is thus, impractical to present the procedure on all the combination of coordinates. We will present the results for the three configuration space coordinates by combining each with its corresponding momentum coordinate. 

\begin{figure*}[!ht] 
	\centering
	\subfigure[]{\includegraphics[width=0.32\textwidth]{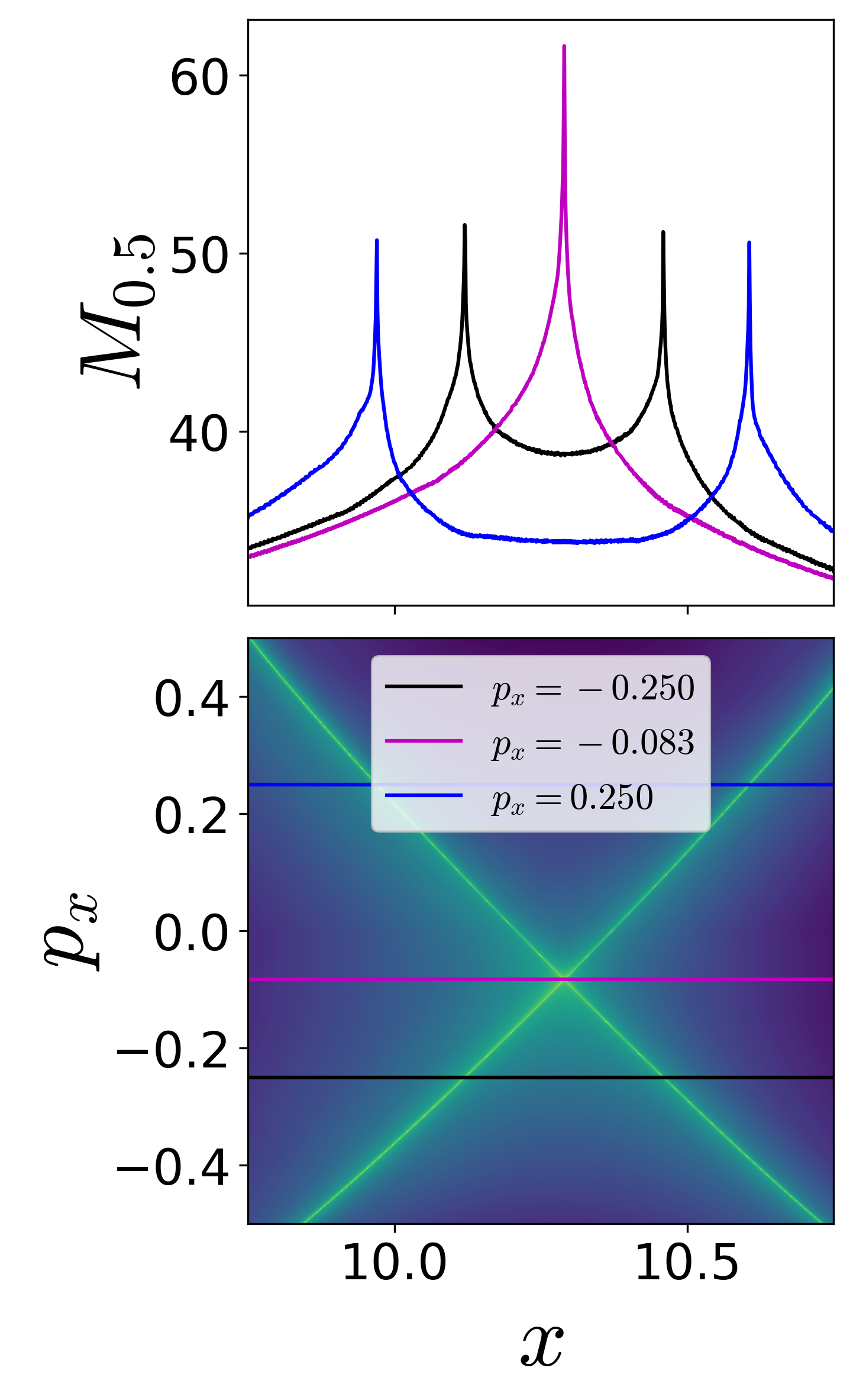}}
	\subfigure[]{\includegraphics[width=0.32\textwidth]{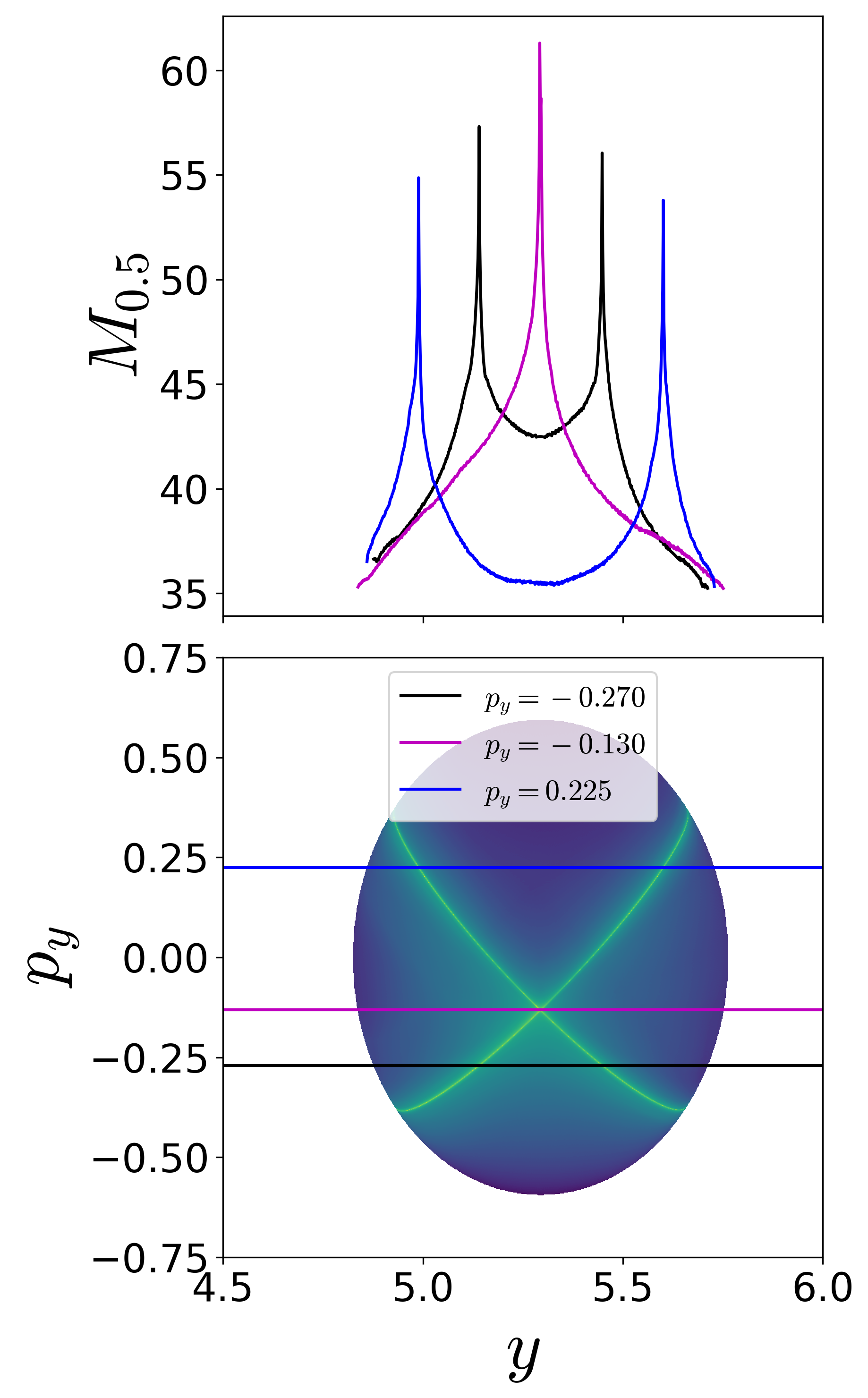}}
	\subfigure[]{\includegraphics[width=0.32\textwidth]{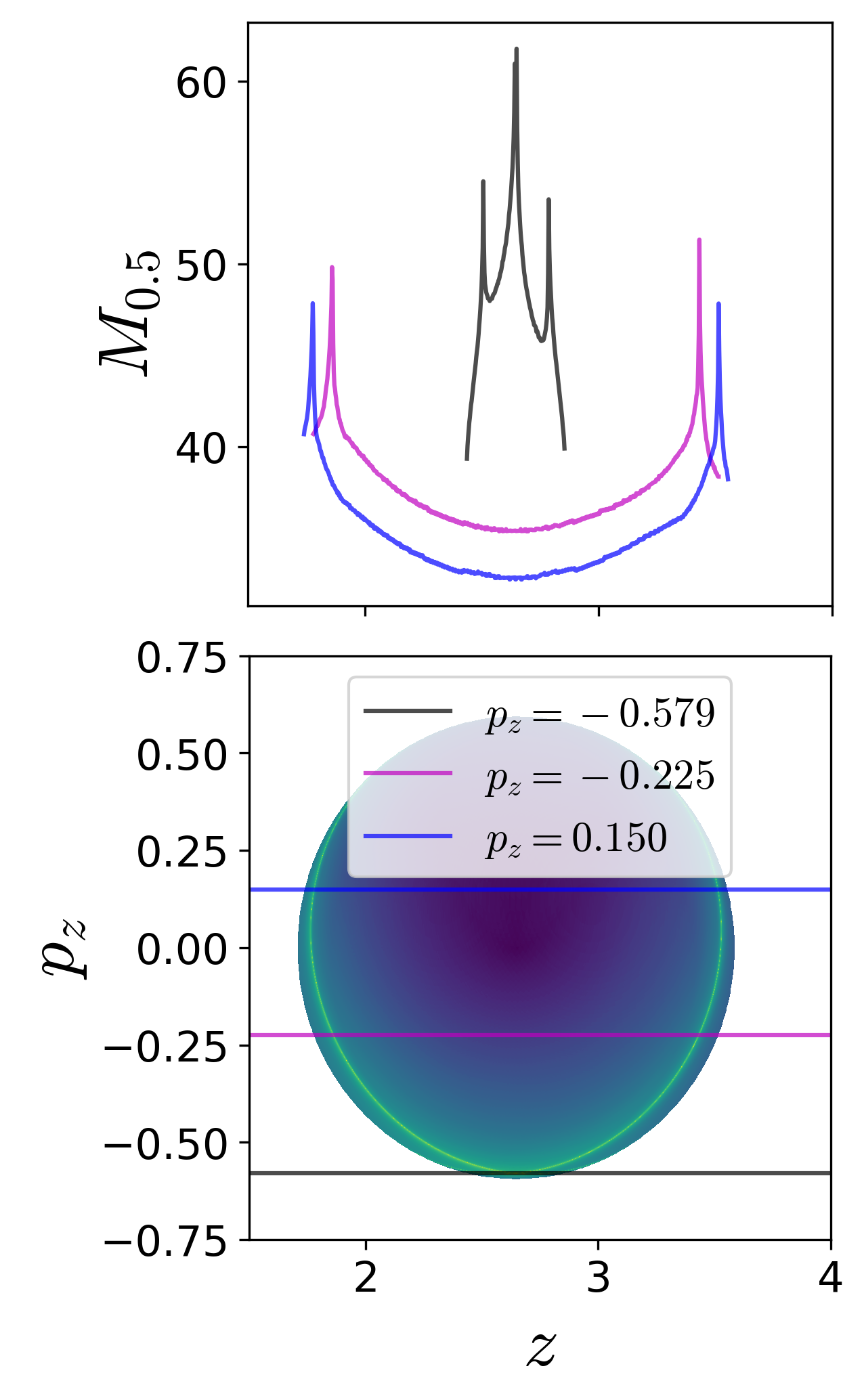}}
	\caption{Detecting points on the NHIM using variable integration time Lagrangian descriptor on the two dimensional surfaces (a) $U_{xp_x}^+$~\eqref{eqn:Barbanis3dof_uxpx}, (b) $U_{yp_y}^+$~\eqref{eqn:Barbanis3dof_uypy}, and (c) $U_{zp_z}^+$~\eqref{eqn:Barbanis3dof_uzpz} at excess energy $\Delta E \approx 0.176$ or total energy $E = 24.000$. For this energy value, the saddle region, as defined in Eqn.~\eqref{eqn:var_time_qs}, is taken to be $q_s = [9,12] \times [2.5,7.5] \times [1,4]$ and $\tau = 50$.}
	\label{fig:Barbanis3dof_M_pxpypz}
\end{figure*}

On these isoenergetic surfaces, we compute the variable integration time Lagrangian descriptor for small excess energy, $\Delta E \approx 0.176$, or total energy $E = 24.000$, and show the contour maps in Fig.~\ref{fig:Barbanis3dof_M_pxpypz}. The maxima identifying the points on the NHIM and its invariant manifolds can be visualized using one dimensional slices for constant momenta. This indicates clearly the initial conditions in the phase space (points on the isoenergetic two dimensional surfaces in $\mathbb{R}^6$, for example~\eqref{eqn:Barbanis3dof_uxpx}) that do not leave the saddle region. 
%\textemdash a simple trajectory diagnostic for detecting points on high dimensional manifolds. 

%and (d-f) $E = 30.000$

\section{\label{sec:summ}Conclusions} %Summary and Outlook

In this article, we discussed a trajectory diagnostic method as a low dimensional probe of high 
dimensional invariant manifolds in 2 and 3 DoF nonlinear Hamiltonian systems. This trajectory 
diagnostic --- Lagrangian descriptor (LD) --- can represent a geometric property of interest in a system with escape/transition and features, that is minima or maxima, in its contour map identify points on the high dimensional invariant manifolds. 

Comparing the points on the NHIM in 2 DoF system obtained using the LD method with differential correction and numerical continuation, we also verified the method for a nonlinear autonomous system following our previous work on decoupled and coupled 2 and 3 DoF linear system~\cite{naik2019finding}. 
The results on 3 DoF system are also congruent with what one expects for an extended problem of the 2 DoF coupled harmonic potential. In addition, the LD method based detection of NHIM is simple to implement and quickly provides a lay of the dynamical land which is a preliminirary step in applying phase space transport to problems in physics and chemistry. 
This method can also be used to set up starting guess for other numerical procedures which rely on 
good initial guess or can also be used in conjunction with machine learning 
methods for rendering the smooth pieces of NHIM~\cite{bardakcioglu2018,feldmaier_invariant_2019}.

\begin{acknowledgments}
We acknowledge the support of EPSRC Grant No. ~EP/P021123/1 and  ONR Grant No.~N00014-01-1-0769. We 
would like to thank Dmitry Zhdanov for stimulating discussions.
\end{acknowledgments}

% \bibliography{LD,Ham_dyn,reaction_dynamics,ross_naik_refs}
\bibliographystyle{aipnum4-1} 
\bibliography{Finding_NHIM_LD_nonlinear}

\providecommand{\noopsort}[1]{}\providecommand{\singleletter}[1]{#1}%
\begin{thebibliography}{74}%
\makeatletter
\providecommand \@ifxundefined [1]{%
 \@ifx{#1\undefined}
}%
\providecommand \@ifnum [1]{%
 \ifnum #1\expandafter \@firstoftwo
 \else \expandafter \@secondoftwo
 \fi
}%
\providecommand \@ifx [1]{%
 \ifx #1\expandafter \@firstoftwo
 \else \expandafter \@secondoftwo
 \fi
}%
\providecommand \natexlab [1]{#1}%
\providecommand \enquote  [1]{``#1''}%
\providecommand \bibnamefont  [1]{#1}%
\providecommand \bibfnamefont [1]{#1}%
\providecommand \citenamefont [1]{#1}%
\providecommand \href@noop [0]{\@secondoftwo}%
\providecommand \href [0]{\begingroup \@sanitize@url \@href}%
\providecommand \@href[1]{\@@startlink{#1}\@@href}%
\providecommand \@@href[1]{\endgroup#1\@@endlink}%
\providecommand \@sanitize@url [0]{\catcode `\\12\catcode `\$12\catcode
  `\&12\catcode `\#12\catcode `\^12\catcode `\_12\catcode `\%12\relax}%
\providecommand \@@startlink[1]{}%
\providecommand \@@endlink[0]{}%
\providecommand \url  [0]{\begingroup\@sanitize@url \@url }%
\providecommand \@url [1]{\endgroup\@href {#1}{\urlprefix }}%
\providecommand \urlprefix  [0]{URL }%
\providecommand \Eprint [0]{\href }%
\providecommand \doibase [0]{http://dx.doi.org/}%
\providecommand \selectlanguage [0]{\@gobble}%
\providecommand \bibinfo  [0]{\@secondoftwo}%
\providecommand \bibfield  [0]{\@secondoftwo}%
\providecommand \translation [1]{[#1]}%
\providecommand \BibitemOpen [0]{}%
\providecommand \bibitemStop [0]{}%
\providecommand \bibitemNoStop [0]{.\EOS\space}%
\providecommand \EOS [0]{\spacefactor3000\relax}%
\providecommand \BibitemShut  [1]{\csname bibitem#1\endcsname}%
\let\auto@bib@innerbib\@empty
%</preamble>
\bibitem [{\citenamefont {Komatsuzaki}\ and\ \citenamefont
  {Berry}(2001)}]{Komatsuzaki2001}%
  \BibitemOpen
  \bibfield  {author} {\bibinfo {author} {\bibfnamefont {T.}~\bibnamefont
  {Komatsuzaki}}\ and\ \bibinfo {author} {\bibfnamefont {R.~S.}\ \bibnamefont
  {Berry}},\ }\href@noop {} {\bibfield  {journal} {\bibinfo  {journal} {P.
  Natl. Acad. Sci. USA}\ }\textbf {\bibinfo {volume} {98}},\ \bibinfo {pages}
  {7666} (\bibinfo {year} {2001})}\BibitemShut {NoStop}%
\bibitem [{\citenamefont {Komatsuzaki}\ and\ \citenamefont
  {Berry}(1999)}]{Komatsuzaki1999}%
  \BibitemOpen
  \bibfield  {author} {\bibinfo {author} {\bibfnamefont {T.}~\bibnamefont
  {Komatsuzaki}}\ and\ \bibinfo {author} {\bibfnamefont {R.~S.}\ \bibnamefont
  {Berry}},\ }\href@noop {} {\bibfield  {journal} {\bibinfo  {journal} {The
  Journal of Chemical Physics}\ }\textbf {\bibinfo {volume} {110}},\ \bibinfo
  {pages} {9160} (\bibinfo {year} {1999})}\BibitemShut {NoStop}%
\bibitem [{\citenamefont {Wiggins}\ \emph {et~al.}(2001)\citenamefont
  {Wiggins}, \citenamefont {Wiesenfeld}, \citenamefont {Jaff\'e},\ and\
  \citenamefont {Uzer}}]{WiWiJaUz2001}%
  \BibitemOpen
  \bibfield  {author} {\bibinfo {author} {\bibfnamefont {S.}~\bibnamefont
  {Wiggins}}, \bibinfo {author} {\bibfnamefont {L.}~\bibnamefont {Wiesenfeld}},
  \bibinfo {author} {\bibfnamefont {C.}~\bibnamefont {Jaff\'e}}, \ and\
  \bibinfo {author} {\bibfnamefont {T.}~\bibnamefont {Uzer}},\ }\href@noop {}
  {\bibfield  {journal} {\bibinfo  {journal} {Phys. Rev. Lett.}\ }\textbf
  {\bibinfo {volume} {86}},\ \bibinfo {pages} {5478} (\bibinfo {year}
  {2001})}\BibitemShut {NoStop}%
\bibitem [{\citenamefont {Jaff\'e}, \citenamefont {Farrelly},\ and\
  \citenamefont {Uzer}(2000)}]{JaFaUz2000}%
  \BibitemOpen
  \bibfield  {author} {\bibinfo {author} {\bibfnamefont {C.}~\bibnamefont
  {Jaff\'e}}, \bibinfo {author} {\bibfnamefont {D.}~\bibnamefont {Farrelly}}, \
  and\ \bibinfo {author} {\bibfnamefont {T.}~\bibnamefont {Uzer}},\ }\href@noop
  {} {\bibfield  {journal} {\bibinfo  {journal} {Phys. Rev. Lett.}\ }\textbf
  {\bibinfo {volume} {84}},\ \bibinfo {pages} {610} (\bibinfo {year}
  {2000})}\BibitemShut {NoStop}%
\bibitem [{\citenamefont {Eckhardt}(1995)}]{Eckhardt1995}%
  \BibitemOpen
  \bibfield  {author} {\bibinfo {author} {\bibfnamefont {B.}~\bibnamefont
  {Eckhardt}},\ }\href@noop {} {\bibfield  {journal} {\bibinfo  {journal} {J.
  Phys. A-Math. Gen.}\ }\textbf {\bibinfo {volume} {28}},\ \bibinfo {pages}
  {3469} (\bibinfo {year} {1995})}\BibitemShut {NoStop}%
\bibitem [{\citenamefont {Collins}, \citenamefont {Ezra},\ and\ \citenamefont
  {Wiggins}(2012)}]{Collins2012}%
  \BibitemOpen
  \bibfield  {author} {\bibinfo {author} {\bibfnamefont {P.}~\bibnamefont
  {Collins}}, \bibinfo {author} {\bibfnamefont {G.~S.}\ \bibnamefont {Ezra}}, \
  and\ \bibinfo {author} {\bibfnamefont {S.}~\bibnamefont {Wiggins}},\
  }\href@noop {} {\bibfield  {journal} {\bibinfo  {journal} {Phys. Rev. E}\
  }\textbf {\bibinfo {volume} {86}},\ \bibinfo {pages} {056218} (\bibinfo
  {year} {2012})}\BibitemShut {NoStop}%
\bibitem [{\citenamefont {Zhong}, \citenamefont {Virgin},\ and\ \citenamefont
  {Ross}(2018)}]{ZhViRo2018}%
  \BibitemOpen
  \bibfield  {author} {\bibinfo {author} {\bibfnamefont {J.}~\bibnamefont
  {Zhong}}, \bibinfo {author} {\bibfnamefont {L.~N.}\ \bibnamefont {Virgin}}, \
  and\ \bibinfo {author} {\bibfnamefont {S.~D.}\ \bibnamefont {Ross}},\
  }\href@noop {} {\bibfield  {journal} {\bibinfo  {journal} {Int. J. Mech.
  Sci.}\ }\textbf {\bibinfo {volume} {000}},\ \bibinfo {pages} {1} (\bibinfo
  {year} {2018})}\BibitemShut {NoStop}%
\bibitem [{\citenamefont {Virgin}(1989)}]{Virgin1989}%
  \BibitemOpen
  \bibfield  {author} {\bibinfo {author} {\bibfnamefont {L.~N.}\ \bibnamefont
  {Virgin}},\ }\href@noop {} {\bibfield  {journal} {\bibinfo  {journal}
  {Dynamics and Stability of Systems}\ }\textbf {\bibinfo {volume} {4}},\
  \bibinfo {pages} {56} (\bibinfo {year} {1989})}\BibitemShut {NoStop}%
\bibitem [{\citenamefont {Thompson}\ and\ \citenamefont {{de
  Souza}}(1996)}]{ThDe1996}%
  \BibitemOpen
  \bibfield  {author} {\bibinfo {author} {\bibfnamefont {J.~M.~T.}\
  \bibnamefont {Thompson}}\ and\ \bibinfo {author} {\bibfnamefont {J.~R.}\
  \bibnamefont {{de Souza}}},\ }\href@noop {} {\bibfield  {journal} {\bibinfo
  {journal} {Proc. R. Soc. Lond. A}\ }\textbf {\bibinfo {volume} {452}},\
  \bibinfo {pages} {2527} (\bibinfo {year} {1996})}\BibitemShut {NoStop}%
\bibitem [{\citenamefont {Naik}\ and\ \citenamefont {Ross}(2017)}]{NaRo2017}%
  \BibitemOpen
  \bibfield  {author} {\bibinfo {author} {\bibfnamefont {S.}~\bibnamefont
  {Naik}}\ and\ \bibinfo {author} {\bibfnamefont {S.~D.}\ \bibnamefont
  {Ross}},\ }\href@noop {} {\bibfield  {journal} {\bibinfo  {journal} {Commun.
  Nonlinear Sci.}\ }\textbf {\bibinfo {volume} {47}},\ \bibinfo {pages} {48}
  (\bibinfo {year} {2017})}\BibitemShut {NoStop}%
\bibitem [{\citenamefont {Jaff\'e}\ \emph {et~al.}(2002)\citenamefont
  {Jaff\'e}, \citenamefont {Ross}, \citenamefont {Lo}, \citenamefont {Marsden},
  \citenamefont {Farrelly},\ and\ \citenamefont {Uzer}}]{JaRoLoMaFaUz2002}%
  \BibitemOpen
  \bibfield  {author} {\bibinfo {author} {\bibfnamefont {C.}~\bibnamefont
  {Jaff\'e}}, \bibinfo {author} {\bibfnamefont {S.~D.}\ \bibnamefont {Ross}},
  \bibinfo {author} {\bibfnamefont {M.~W.}\ \bibnamefont {Lo}}, \bibinfo
  {author} {\bibfnamefont {J.~E.}\ \bibnamefont {Marsden}}, \bibinfo {author}
  {\bibfnamefont {D.}~\bibnamefont {Farrelly}}, \ and\ \bibinfo {author}
  {\bibfnamefont {T.}~\bibnamefont {Uzer}},\ }\href@noop {} {\bibfield
  {journal} {\bibinfo  {journal} {Physical Review Letters}\ }\textbf {\bibinfo
  {volume} {89}},\ \bibinfo {pages} {011101} (\bibinfo {year}
  {2002})}\BibitemShut {NoStop}%
\bibitem [{\citenamefont {Dellnitz}\ \emph {et~al.}(2005)\citenamefont
  {Dellnitz}, \citenamefont {Junge}, \citenamefont {Lo}, \citenamefont
  {Marsden}, \citenamefont {Padberg}, \citenamefont {Preis}, \citenamefont
  {Ross},\ and\ \citenamefont {Thiere}}]{DeJuLoMaPaPrRoTh2005}%
  \BibitemOpen
  \bibfield  {author} {\bibinfo {author} {\bibfnamefont {M.}~\bibnamefont
  {Dellnitz}}, \bibinfo {author} {\bibfnamefont {O.}~\bibnamefont {Junge}},
  \bibinfo {author} {\bibfnamefont {M.~W.}\ \bibnamefont {Lo}}, \bibinfo
  {author} {\bibfnamefont {J.~E.}\ \bibnamefont {Marsden}}, \bibinfo {author}
  {\bibfnamefont {K.}~\bibnamefont {Padberg}}, \bibinfo {author} {\bibfnamefont
  {R.}~\bibnamefont {Preis}}, \bibinfo {author} {\bibfnamefont {S.~D.}\
  \bibnamefont {Ross}}, \ and\ \bibinfo {author} {\bibfnamefont
  {B.}~\bibnamefont {Thiere}},\ }\href@noop {} {\bibfield  {journal} {\bibinfo
  {journal} {Physical Review Letters}\ }\textbf {\bibinfo {volume} {94}},\
  \bibinfo {pages} {231102} (\bibinfo {year} {2005})}\BibitemShut {NoStop}%
\bibitem [{\citenamefont {Ross}(2003)}]{Ross2003}%
  \BibitemOpen
  \bibfield  {author} {\bibinfo {author} {\bibfnamefont {S.~D.}\ \bibnamefont
  {Ross}},\ }in\ \href@noop {} {\emph {\bibinfo {booktitle} {Libration Point
  Orbits and Applications}}},\ \bibinfo {editor} {edited by\ \bibinfo {editor}
  {\bibfnamefont {G.}~\bibnamefont {G\'omez}}, \bibinfo {editor} {\bibfnamefont
  {M.~W.}\ \bibnamefont {Lo}}, \ and\ \bibinfo {editor} {\bibfnamefont {J.~J.}\
  \bibnamefont {Masdemont}}}\ (\bibinfo  {publisher} {World Scientific},\
  \bibinfo {year} {2003})\ pp.\ \bibinfo {pages} {637--652}\BibitemShut
  {NoStop}%
\bibitem [{\citenamefont {de~Oliveira}\ \emph {et~al.}(2002)\citenamefont
  {de~Oliveira}, \citenamefont {{Ozorio de Almeida}}, \citenamefont
  {{Dami{\~{a}}o Soares}},\ and\ \citenamefont {Tonini}}]{DeOliveira2002}%
  \BibitemOpen
  \bibfield  {author} {\bibinfo {author} {\bibfnamefont {H.~P.}\ \bibnamefont
  {de~Oliveira}}, \bibinfo {author} {\bibfnamefont {A.~M.}\ \bibnamefont
  {{Ozorio de Almeida}}}, \bibinfo {author} {\bibfnamefont {I.}~\bibnamefont
  {{Dami{\~{a}}o Soares}}}, \ and\ \bibinfo {author} {\bibfnamefont {E.~V.}\
  \bibnamefont {Tonini}},\ }\href@noop {} {\bibfield  {journal} {\bibinfo
  {journal} {Phys. Rev. D}\ }\textbf {\bibinfo {volume} {65}},\ \bibinfo
  {pages} {9} (\bibinfo {year} {2002})}\BibitemShut {NoStop}%
\bibitem [{\citenamefont {Patra}\ and\ \citenamefont
  {Keshavamurthy}(2015)}]{patra_classical-quantum_2015}%
  \BibitemOpen
  \bibfield  {author} {\bibinfo {author} {\bibfnamefont {S.}~\bibnamefont
  {Patra}}\ and\ \bibinfo {author} {\bibfnamefont {S.}~\bibnamefont
  {Keshavamurthy}},\ }\href@noop {} {\bibfield  {journal} {\bibinfo  {journal}
  {Chemical Physics Letters}\ }\textbf {\bibinfo {volume} {634}},\ \bibinfo
  {pages} {1} (\bibinfo {year} {2015})}\BibitemShut {NoStop}%
\bibitem [{\citenamefont {Patra}\ and\ \citenamefont
  {Keshavamurthy}(2018)}]{patra_detecting_2018}%
  \BibitemOpen
  \bibfield  {author} {\bibinfo {author} {\bibfnamefont {S.}~\bibnamefont
  {Patra}}\ and\ \bibinfo {author} {\bibfnamefont {S.}~\bibnamefont
  {Keshavamurthy}},\ }\href@noop {} {\bibfield  {journal} {\bibinfo  {journal}
  {Physical Chemistry Chemical Physics}\ }\textbf {\bibinfo {volume} {20}},\
  \bibinfo {pages} {4970} (\bibinfo {year} {2018})}\BibitemShut {NoStop}%
\bibitem [{\citenamefont {Naik}, \citenamefont {Garc{\'\i}a-Garrido},\ and\
  \citenamefont {Wiggins}(2019)}]{naik2019finding}%
  \BibitemOpen
  \bibfield  {author} {\bibinfo {author} {\bibfnamefont {S.}~\bibnamefont
  {Naik}}, \bibinfo {author} {\bibfnamefont {V.~J.}\ \bibnamefont
  {Garc{\'\i}a-Garrido}}, \ and\ \bibinfo {author} {\bibfnamefont
  {S.}~\bibnamefont {Wiggins}},\ }\href@noop {} {\bibfield  {journal} {\bibinfo
   {journal} {arXiv preprint arXiv:1903.10264}\ } (\bibinfo {year}
  {2019})}\BibitemShut {NoStop}%
\bibitem [{\citenamefont {Barbanis}(1966)}]{barbanis_isolating_1966}%
  \BibitemOpen
  \bibfield  {author} {\bibinfo {author} {\bibfnamefont {B.}~\bibnamefont
  {Barbanis}},\ }\href@noop {} {\bibfield  {journal} {\bibinfo  {journal} {The
  Astronomical Journal}\ }\textbf {\bibinfo {volume} {71}},\ \bibinfo {pages}
  {415} (\bibinfo {year} {1966})}\BibitemShut {NoStop}%
\bibitem [{\citenamefont {Brumer}\ and\ \citenamefont
  {Duff}(1976)}]{brumer_variational_1976}%
  \BibitemOpen
  \bibfield  {author} {\bibinfo {author} {\bibfnamefont {P.}~\bibnamefont
  {Brumer}}\ and\ \bibinfo {author} {\bibfnamefont {J.~W.}\ \bibnamefont
  {Duff}},\ }\href@noop {} {\bibfield  {journal} {\bibinfo  {journal} {The
  Journal of Chemical Physics}\ }\textbf {\bibinfo {volume} {65}},\ \bibinfo
  {pages} {3566} (\bibinfo {year} {1976})}\BibitemShut {NoStop}%
\bibitem [{\citenamefont {Davis}\ and\ \citenamefont
  {Heller}(1979)}]{davis_semiclassical_1979}%
  \BibitemOpen
  \bibfield  {author} {\bibinfo {author} {\bibfnamefont {M.~J.}\ \bibnamefont
  {Davis}}\ and\ \bibinfo {author} {\bibfnamefont {E.~J.}\ \bibnamefont
  {Heller}},\ }\href@noop {} {\bibfield  {journal} {\bibinfo  {journal} {The
  Journal of Chemical Physics}\ }\textbf {\bibinfo {volume} {71}},\ \bibinfo
  {pages} {3383} (\bibinfo {year} {1979})}\BibitemShut {NoStop}%
\bibitem [{\citenamefont {Heller}, \citenamefont {Stechel},\ and\ \citenamefont
  {Davis}(1980{\natexlab{a}})}]{heller_molecular_1980}%
  \BibitemOpen
  \bibfield  {author} {\bibinfo {author} {\bibfnamefont {E.~J.}\ \bibnamefont
  {Heller}}, \bibinfo {author} {\bibfnamefont {E.~B.}\ \bibnamefont {Stechel}},
  \ and\ \bibinfo {author} {\bibfnamefont {M.~J.}\ \bibnamefont {Davis}},\
  }\href@noop {} {\bibfield  {journal} {\bibinfo  {journal} {The Journal of
  Chemical Physics}\ }\textbf {\bibinfo {volume} {73}},\ \bibinfo {pages}
  {4720} (\bibinfo {year} {1980}{\natexlab{a}})}\BibitemShut {NoStop}%
\bibitem [{\citenamefont {Waite}\ and\ \citenamefont
  {Miller}(1981)}]{waite_mode_1981}%
  \BibitemOpen
  \bibfield  {author} {\bibinfo {author} {\bibfnamefont {B.~A.}\ \bibnamefont
  {Waite}}\ and\ \bibinfo {author} {\bibfnamefont {W.~H.}\ \bibnamefont
  {Miller}},\ }\href@noop {} {\bibfield  {journal} {\bibinfo  {journal} {The
  Journal of Chemical Physics}\ }\textbf {\bibinfo {volume} {74}},\ \bibinfo
  {pages} {3910} (\bibinfo {year} {1981})}\BibitemShut {NoStop}%
\bibitem [{\citenamefont {Kosloff}\ and\ \citenamefont
  {Rice}(1981)}]{kosloff_dynamical_1981}%
  \BibitemOpen
  \bibfield  {author} {\bibinfo {author} {\bibfnamefont {R.}~\bibnamefont
  {Kosloff}}\ and\ \bibinfo {author} {\bibfnamefont {S.~A.}\ \bibnamefont
  {Rice}},\ }\href@noop {} {\bibfield  {journal} {\bibinfo  {journal} {The
  Journal of Chemical Physics}\ }\textbf {\bibinfo {volume} {74}},\ \bibinfo
  {pages} {1947} (\bibinfo {year} {1981})}\BibitemShut {NoStop}%
\bibitem [{\citenamefont {Contopoulos}\ and\ \citenamefont
  {Magnenat}(1985)}]{contopoulos_simple_1985}%
  \BibitemOpen
  \bibfield  {author} {\bibinfo {author} {\bibfnamefont {G.}~\bibnamefont
  {Contopoulos}}\ and\ \bibinfo {author} {\bibfnamefont {P.}~\bibnamefont
  {Magnenat}},\ }\href@noop {} {\bibfield  {journal} {\bibinfo  {journal}
  {Celestial Mechanics}\ }\textbf {\bibinfo {volume} {37}},\ \bibinfo {pages}
  {387} (\bibinfo {year} {1985})}\BibitemShut {NoStop}%
\bibitem [{\citenamefont {Founargiotakis}\ \emph {et~al.}(1989)\citenamefont
  {Founargiotakis}, \citenamefont {Farantos}, \citenamefont {Contopoulos},\
  and\ \citenamefont {Polymilis}}]{founargiotakis_periodic_1989}%
  \BibitemOpen
  \bibfield  {author} {\bibinfo {author} {\bibfnamefont {M.}~\bibnamefont
  {Founargiotakis}}, \bibinfo {author} {\bibfnamefont {S.~C.}\ \bibnamefont
  {Farantos}}, \bibinfo {author} {\bibfnamefont {G.}~\bibnamefont
  {Contopoulos}}, \ and\ \bibinfo {author} {\bibfnamefont {C.}~\bibnamefont
  {Polymilis}},\ }\href@noop {} {\bibfield  {journal} {\bibinfo  {journal} {The
  Journal of Chemical Physics}\ }\textbf {\bibinfo {volume} {91}},\ \bibinfo
  {pages} {1389} (\bibinfo {year} {1989})}\BibitemShut {NoStop}%
\bibitem [{\citenamefont {Barbanis}(1990)}]{barbanis_escape_1990}%
  \BibitemOpen
  \bibfield  {author} {\bibinfo {author} {\bibfnamefont {B.}~\bibnamefont
  {Barbanis}},\ }\href@noop {} {\bibfield  {journal} {\bibinfo  {journal}
  {Celestial Mechanics and Dynamical Astronomy}\ }\textbf {\bibinfo {volume}
  {48}},\ \bibinfo {pages} {57} (\bibinfo {year} {1990})}\BibitemShut {NoStop}%
\bibitem [{\citenamefont {Babyuk}, \citenamefont {Wyatt},\ and\ \citenamefont
  {Frederick}(2003)}]{babyuk_hydrodynamic_2003}%
  \BibitemOpen
  \bibfield  {author} {\bibinfo {author} {\bibfnamefont {D.}~\bibnamefont
  {Babyuk}}, \bibinfo {author} {\bibfnamefont {R.~E.}\ \bibnamefont {Wyatt}}, \
  and\ \bibinfo {author} {\bibfnamefont {J.~H.}\ \bibnamefont {Frederick}},\
  }\href@noop {} {\bibfield  {journal} {\bibinfo  {journal} {The Journal of
  Chemical Physics}\ }\textbf {\bibinfo {volume} {119}},\ \bibinfo {pages}
  {6482} (\bibinfo {year} {2003})}\BibitemShut {NoStop}%
\bibitem [{\citenamefont {Mitchell}\ \emph
  {et~al.}(2003{\natexlab{a}})\citenamefont {Mitchell}, \citenamefont
  {Handley}, \citenamefont {Tighe}, \citenamefont {Delos},\ and\ \citenamefont
  {Knudson}}]{mitchell_geometry_2003_I}%
  \BibitemOpen
  \bibfield  {author} {\bibinfo {author} {\bibfnamefont {K.~A.}\ \bibnamefont
  {Mitchell}}, \bibinfo {author} {\bibfnamefont {J.~P.}\ \bibnamefont
  {Handley}}, \bibinfo {author} {\bibfnamefont {B.}~\bibnamefont {Tighe}},
  \bibinfo {author} {\bibfnamefont {J.~B.}\ \bibnamefont {Delos}}, \ and\
  \bibinfo {author} {\bibfnamefont {S.~K.}\ \bibnamefont {Knudson}},\
  }\href@noop {} {\bibfield  {journal} {\bibinfo  {journal} {Chaos: An
  Interdisciplinary Journal of Nonlinear Science}\ }\textbf {\bibinfo {volume}
  {13}},\ \bibinfo {pages} {880} (\bibinfo {year}
  {2003}{\natexlab{a}})}\BibitemShut {NoStop}%
\bibitem [{\citenamefont {Mitchell}\ \emph
  {et~al.}(2003{\natexlab{b}})\citenamefont {Mitchell}, \citenamefont
  {Handley}, \citenamefont {Delos},\ and\ \citenamefont
  {Knudson}}]{mitchell_geometry_2003_II}%
  \BibitemOpen
  \bibfield  {author} {\bibinfo {author} {\bibfnamefont {K.~A.}\ \bibnamefont
  {Mitchell}}, \bibinfo {author} {\bibfnamefont {J.~P.}\ \bibnamefont
  {Handley}}, \bibinfo {author} {\bibfnamefont {J.~B.}\ \bibnamefont {Delos}},
  \ and\ \bibinfo {author} {\bibfnamefont {S.~K.}\ \bibnamefont {Knudson}},\
  }\href@noop {} {\bibfield  {journal} {\bibinfo  {journal} {Chaos: An
  Interdisciplinary Journal of Nonlinear Science}\ }\textbf {\bibinfo {volume}
  {13}},\ \bibinfo {pages} {892} (\bibinfo {year}
  {2003}{\natexlab{b}})}\BibitemShut {NoStop}%
\bibitem [{\citenamefont {Mitchell}\ \emph
  {et~al.}(2004{\natexlab{a}})\citenamefont {Mitchell}, \citenamefont
  {Handley}, \citenamefont {Tighe}, \citenamefont {Flower},\ and\ \citenamefont
  {Delos}}]{mitchell_chaos-induced_2004}%
  \BibitemOpen
  \bibfield  {author} {\bibinfo {author} {\bibfnamefont {K.~A.}\ \bibnamefont
  {Mitchell}}, \bibinfo {author} {\bibfnamefont {J.~P.}\ \bibnamefont
  {Handley}}, \bibinfo {author} {\bibfnamefont {B.}~\bibnamefont {Tighe}},
  \bibinfo {author} {\bibfnamefont {A.}~\bibnamefont {Flower}}, \ and\ \bibinfo
  {author} {\bibfnamefont {J.~B.}\ \bibnamefont {Delos}},\ }\href@noop {}
  {\bibfield  {journal} {\bibinfo  {journal} {Physical Review Letters}\
  }\textbf {\bibinfo {volume} {92}} (\bibinfo {year}
  {2004}{\natexlab{a}})}\BibitemShut {NoStop}%
\bibitem [{\citenamefont {De~Leon}\ and\ \citenamefont
  {Berne}(1981)}]{de_leon_intramolecular_1981}%
  \BibitemOpen
  \bibfield  {author} {\bibinfo {author} {\bibfnamefont {N.}~\bibnamefont
  {De~Leon}}\ and\ \bibinfo {author} {\bibfnamefont {B.~J.}\ \bibnamefont
  {Berne}},\ }\href@noop {} {\bibfield  {journal} {\bibinfo  {journal} {The
  Journal of Chemical Physics}\ }\textbf {\bibinfo {volume} {75}},\ \bibinfo
  {pages} {3495} (\bibinfo {year} {1981})}\BibitemShut {NoStop}%
\bibitem [{\citenamefont {Mitchell}\ \emph
  {et~al.}(2004{\natexlab{b}})\citenamefont {Mitchell}, \citenamefont
  {Handley}, \citenamefont {Tighe}, \citenamefont {Flower},\ and\ \citenamefont
  {Delos}}]{mitchell_analysis_2004}%
  \BibitemOpen
  \bibfield  {author} {\bibinfo {author} {\bibfnamefont {K.~A.}\ \bibnamefont
  {Mitchell}}, \bibinfo {author} {\bibfnamefont {J.~P.}\ \bibnamefont
  {Handley}}, \bibinfo {author} {\bibfnamefont {B.}~\bibnamefont {Tighe}},
  \bibinfo {author} {\bibfnamefont {A.}~\bibnamefont {Flower}}, \ and\ \bibinfo
  {author} {\bibfnamefont {J.~B.}\ \bibnamefont {Delos}},\ }\href@noop {}
  {\bibfield  {journal} {\bibinfo  {journal} {Physical Review A}\ }\textbf
  {\bibinfo {volume} {70}} (\bibinfo {year} {2004}{\natexlab{b}})}\BibitemShut
  {NoStop}%
\bibitem [{\citenamefont {Mitchell}\ and\ \citenamefont
  {Ilan}(2009)}]{mitchell_nonlinear_2009}%
  \BibitemOpen
  \bibfield  {author} {\bibinfo {author} {\bibfnamefont {K.~A.}\ \bibnamefont
  {Mitchell}}\ and\ \bibinfo {author} {\bibfnamefont {B.}~\bibnamefont
  {Ilan}},\ }\href@noop {} {\bibfield  {journal} {\bibinfo  {journal} {Physical
  Review A}\ }\textbf {\bibinfo {volume} {80}} (\bibinfo {year}
  {2009})}\BibitemShut {NoStop}%
\bibitem [{\citenamefont {Mitchell}\ and\ \citenamefont
  {Delos}(2007)}]{mitchell_structure_2007}%
  \BibitemOpen
  \bibfield  {author} {\bibinfo {author} {\bibfnamefont {K.~A.}\ \bibnamefont
  {Mitchell}}\ and\ \bibinfo {author} {\bibfnamefont {J.~B.}\ \bibnamefont
  {Delos}},\ }\href@noop {} {\bibfield  {journal} {\bibinfo  {journal} {Physica
  D: Nonlinear Phenomena}\ }\textbf {\bibinfo {volume} {229}},\ \bibinfo
  {pages} {9} (\bibinfo {year} {2007})}\BibitemShut {NoStop}%
\bibitem [{\citenamefont {Wang}\ \emph {et~al.}(2010)\citenamefont {Wang},
  \citenamefont {Yang}, \citenamefont {Liu}, \citenamefont {Liu}, \citenamefont
  {Zhan},\ and\ \citenamefont {Delos}}]{wang_photoionization_2010}%
  \BibitemOpen
  \bibfield  {author} {\bibinfo {author} {\bibfnamefont {L.}~\bibnamefont
  {Wang}}, \bibinfo {author} {\bibfnamefont {H.~F.}\ \bibnamefont {Yang}},
  \bibinfo {author} {\bibfnamefont {X.~J.}\ \bibnamefont {Liu}}, \bibinfo
  {author} {\bibfnamefont {H.~P.}\ \bibnamefont {Liu}}, \bibinfo {author}
  {\bibfnamefont {M.~S.}\ \bibnamefont {Zhan}}, \ and\ \bibinfo {author}
  {\bibfnamefont {J.~B.}\ \bibnamefont {Delos}},\ }\href@noop {} {\bibfield
  {journal} {\bibinfo  {journal} {Physical Review A}\ }\textbf {\bibinfo
  {volume} {82}} (\bibinfo {year} {2010})}\BibitemShut {NoStop}%
\bibitem [{\citenamefont {Madrid}\ and\ \citenamefont
  {Mancho}(2009)}]{madrid2009}%
  \BibitemOpen
  \bibfield  {author} {\bibinfo {author} {\bibfnamefont {J.~A.~J.}\
  \bibnamefont {Madrid}}\ and\ \bibinfo {author} {\bibfnamefont {A.~M.}\
  \bibnamefont {Mancho}},\ }\href@noop {} {\bibfield  {journal} {\bibinfo
  {journal} {Chaos}\ }\textbf {\bibinfo {volume} {19}},\ \bibinfo {pages}
  {013111} (\bibinfo {year} {2009})}\BibitemShut {NoStop}%
\bibitem [{\citenamefont {Mendoza}\ and\ \citenamefont
  {Mancho}(2010)}]{mendoza2010}%
  \BibitemOpen
  \bibfield  {author} {\bibinfo {author} {\bibfnamefont {C.}~\bibnamefont
  {Mendoza}}\ and\ \bibinfo {author} {\bibfnamefont {A.~M.}\ \bibnamefont
  {Mancho}},\ }\href@noop {} {\bibfield  {journal} {\bibinfo  {journal} {Phys.
  Rev. Lett.}\ }\textbf {\bibinfo {volume} {105}},\ \bibinfo {pages} {038501}
  (\bibinfo {year} {2010})}\BibitemShut {NoStop}%
\bibitem [{\citenamefont {Mancho}\ \emph {et~al.}(2013)\citenamefont {Mancho},
  \citenamefont {Wiggins}, \citenamefont {Curbelo},\ and\ \citenamefont
  {Mendoza}}]{mancho2013}%
  \BibitemOpen
  \bibfield  {author} {\bibinfo {author} {\bibfnamefont {A.~M.}\ \bibnamefont
  {Mancho}}, \bibinfo {author} {\bibfnamefont {S.}~\bibnamefont {Wiggins}},
  \bibinfo {author} {\bibfnamefont {J.}~\bibnamefont {Curbelo}}, \ and\
  \bibinfo {author} {\bibfnamefont {C.}~\bibnamefont {Mendoza}},\ }\href@noop
  {} {\bibfield  {journal} {\bibinfo  {journal} {Communications in Nonlinear
  Science and Numerical}\ }\textbf {\bibinfo {volume} {18}},\ \bibinfo {pages}
  {3530} (\bibinfo {year} {2013})}\BibitemShut {NoStop}%
\bibitem [{\citenamefont {Lopesino}\ \emph {et~al.}(2017)\citenamefont
  {Lopesino}, \citenamefont {Balibrea-Iniesta}, \citenamefont
  {Garc\'ia-Garrido}, \citenamefont {Wiggins},\ and\ \citenamefont
  {Mancho}}]{lopesino2017}%
  \BibitemOpen
  \bibfield  {author} {\bibinfo {author} {\bibfnamefont {C.}~\bibnamefont
  {Lopesino}}, \bibinfo {author} {\bibfnamefont {F.}~\bibnamefont
  {Balibrea-Iniesta}}, \bibinfo {author} {\bibfnamefont {V.~J.}\ \bibnamefont
  {Garc\'ia-Garrido}}, \bibinfo {author} {\bibfnamefont {S.}~\bibnamefont
  {Wiggins}}, \ and\ \bibinfo {author} {\bibfnamefont {A.~M.}\ \bibnamefont
  {Mancho}},\ }\href@noop {} {\bibfield  {journal} {\bibinfo  {journal}
  {International Journal of Bifurcation and Chaos}\ }\textbf {\bibinfo {volume}
  {27}},\ \bibinfo {pages} {1730001} (\bibinfo {year} {2017})}\BibitemShut
  {NoStop}%
\bibitem [{\citenamefont {Balibrea-Iniesta}\ \emph {et~al.}(2016)\citenamefont
  {Balibrea-Iniesta}, \citenamefont {Lopesino}, \citenamefont {Wiggins},\ and\
  \citenamefont {Mancho}}]{balibrea2016lagrangian}%
  \BibitemOpen
  \bibfield  {author} {\bibinfo {author} {\bibfnamefont {F.}~\bibnamefont
  {Balibrea-Iniesta}}, \bibinfo {author} {\bibfnamefont {C.}~\bibnamefont
  {Lopesino}}, \bibinfo {author} {\bibfnamefont {S.}~\bibnamefont {Wiggins}}, \
  and\ \bibinfo {author} {\bibfnamefont {A.~M.}\ \bibnamefont {Mancho}},\
  }\href@noop {} {\bibfield  {journal} {\bibinfo  {journal} {International
  Journal of Bifurcation and Chaos}\ }\textbf {\bibinfo {volume} {26}},\
  \bibinfo {pages} {1630036} (\bibinfo {year} {2016})}\BibitemShut {NoStop}%
\bibitem [{\citenamefont {Craven}, \citenamefont {Junginger},\ and\
  \citenamefont {Hernandez}(2017)}]{craven2017lagrangian}%
  \BibitemOpen
  \bibfield  {author} {\bibinfo {author} {\bibfnamefont {G.~T.}\ \bibnamefont
  {Craven}}, \bibinfo {author} {\bibfnamefont {A.}~\bibnamefont {Junginger}}, \
  and\ \bibinfo {author} {\bibfnamefont {R.}~\bibnamefont {Hernandez}},\
  }\href@noop {} {\bibfield  {journal} {\bibinfo  {journal} {Physical Review
  E}\ }\textbf {\bibinfo {volume} {96}},\ \bibinfo {pages} {022222} (\bibinfo
  {year} {2017})}\BibitemShut {NoStop}%
\bibitem [{\citenamefont {Junginger}\ and\ \citenamefont
  {Hernandez}(2016{\natexlab{a}})}]{junginger2016lagrangian}%
  \BibitemOpen
  \bibfield  {author} {\bibinfo {author} {\bibfnamefont {A.}~\bibnamefont
  {Junginger}}\ and\ \bibinfo {author} {\bibfnamefont {R.}~\bibnamefont
  {Hernandez}},\ }\href@noop {} {\bibfield  {journal} {\bibinfo  {journal}
  {Physical Chemistry Chemical Physics}\ }\textbf {\bibinfo {volume} {18}},\
  \bibinfo {pages} {30282} (\bibinfo {year} {2016}{\natexlab{a}})}\BibitemShut
  {NoStop}%
\bibitem [{\citenamefont {de~la C{\'a}mara}\ \emph {et~al.}(2012)\citenamefont
  {de~la C{\'a}mara}, \citenamefont {Mancho}, \citenamefont {Ide},
  \citenamefont {Serrano},\ and\ \citenamefont {Mechoso}}]{amism11}%
  \BibitemOpen
  \bibfield  {author} {\bibinfo {author} {\bibfnamefont {A.}~\bibnamefont
  {de~la C{\'a}mara}}, \bibinfo {author} {\bibfnamefont {A.~M.}\ \bibnamefont
  {Mancho}}, \bibinfo {author} {\bibfnamefont {K.}~\bibnamefont {Ide}},
  \bibinfo {author} {\bibfnamefont {E.}~\bibnamefont {Serrano}}, \ and\
  \bibinfo {author} {\bibfnamefont {C.}~\bibnamefont {Mechoso}},\ }\href@noop
  {} {\bibfield  {journal} {\bibinfo  {journal} {J. Atmos. Sci.}\ }\textbf
  {\bibinfo {volume} {69}},\ \bibinfo {pages} {753} (\bibinfo {year}
  {2012})}\BibitemShut {NoStop}%
\bibitem [{\citenamefont {Mendoza}, \citenamefont {Mancho},\ and\ \citenamefont
  {Wiggins}(2014)}]{mendoza2014}%
  \BibitemOpen
  \bibfield  {author} {\bibinfo {author} {\bibfnamefont {C.}~\bibnamefont
  {Mendoza}}, \bibinfo {author} {\bibfnamefont {A.~M.}\ \bibnamefont {Mancho}},
  \ and\ \bibinfo {author} {\bibfnamefont {S.}~\bibnamefont {Wiggins}},\
  }\href@noop {} {\bibfield  {journal} {\bibinfo  {journal} {Nonlinear
  Processes in Geophysics}\ }\textbf {\bibinfo {volume} {21}},\ \bibinfo
  {pages} {677} (\bibinfo {year} {2014})}\BibitemShut {NoStop}%
\bibitem [{\citenamefont {Garc{\'\i}a-Garrido}, \citenamefont {Mancho},\ and\
  \citenamefont {Wiggins}(2015)}]{ggmwm15}%
  \BibitemOpen
  \bibfield  {author} {\bibinfo {author} {\bibfnamefont {V.~J.}\ \bibnamefont
  {Garc{\'\i}a-Garrido}}, \bibinfo {author} {\bibfnamefont {A.~M.}\
  \bibnamefont {Mancho}}, \ and\ \bibinfo {author} {\bibfnamefont
  {S.}~\bibnamefont {Wiggins}},\ }\href@noop {} {\bibfield  {journal} {\bibinfo
   {journal} {Nonlin. Proc. Geophys.}\ }\textbf {\bibinfo {volume} {22}},\
  \bibinfo {pages} {701} (\bibinfo {year} {2015})}\BibitemShut {NoStop}%
\bibitem [{\citenamefont {Ramos}\ \emph {et~al.}(2018)\citenamefont {Ramos},
  \citenamefont {Garc{\'\i}a-Garrido}, \citenamefont {Mancho}, \citenamefont
  {Wiggins}, \citenamefont {Coca}, \citenamefont {Glenn}, \citenamefont
  {Schofield}, \citenamefont {Kohut}, \citenamefont {Aragon}, \citenamefont
  {Kerfoot}, \citenamefont {Haskins}, \citenamefont {Miles}, \citenamefont
  {Haldeman}, \citenamefont {Strandskov}, \citenamefont {Allsup}, \citenamefont
  {Jones},\ and\ \citenamefont {Shapiro.}}]{ramos2018}%
  \BibitemOpen
  \bibfield  {author} {\bibinfo {author} {\bibfnamefont {A.~G.}\ \bibnamefont
  {Ramos}}, \bibinfo {author} {\bibfnamefont {V.~J.}\ \bibnamefont
  {Garc{\'\i}a-Garrido}}, \bibinfo {author} {\bibfnamefont {A.~M.}\
  \bibnamefont {Mancho}}, \bibinfo {author} {\bibfnamefont {S.}~\bibnamefont
  {Wiggins}}, \bibinfo {author} {\bibfnamefont {J.}~\bibnamefont {Coca}},
  \bibinfo {author} {\bibfnamefont {S.}~\bibnamefont {Glenn}}, \bibinfo
  {author} {\bibfnamefont {O.}~\bibnamefont {Schofield}}, \bibinfo {author}
  {\bibfnamefont {J.}~\bibnamefont {Kohut}}, \bibinfo {author} {\bibfnamefont
  {D.}~\bibnamefont {Aragon}}, \bibinfo {author} {\bibfnamefont
  {J.}~\bibnamefont {Kerfoot}}, \bibinfo {author} {\bibfnamefont
  {T.}~\bibnamefont {Haskins}}, \bibinfo {author} {\bibfnamefont
  {T.}~\bibnamefont {Miles}}, \bibinfo {author} {\bibfnamefont
  {C.}~\bibnamefont {Haldeman}}, \bibinfo {author} {\bibfnamefont
  {N.}~\bibnamefont {Strandskov}}, \bibinfo {author} {\bibfnamefont
  {B.}~\bibnamefont {Allsup}}, \bibinfo {author} {\bibfnamefont
  {C.}~\bibnamefont {Jones}}, \ and\ \bibinfo {author} {\bibfnamefont
  {J.}~\bibnamefont {Shapiro.}},\ }\href@noop {} {\bibfield  {journal}
  {\bibinfo  {journal} {Scientfic Reports}\ }\textbf {\bibinfo {volume} {4}},\
  \bibinfo {pages} {4575} (\bibinfo {year} {2018})}\BibitemShut {NoStop}%
\bibitem [{\citenamefont {Junginger}\ \emph {et~al.}(2016)\citenamefont
  {Junginger}, \citenamefont {Craven}, \citenamefont {Bartsch}, \citenamefont
  {Revuelta}, \citenamefont {Borondo}, \citenamefont {Benito},\ and\
  \citenamefont {Hernandez}}]{junginger2016transition}%
  \BibitemOpen
  \bibfield  {author} {\bibinfo {author} {\bibfnamefont {A.}~\bibnamefont
  {Junginger}}, \bibinfo {author} {\bibfnamefont {G.~T.}\ \bibnamefont
  {Craven}}, \bibinfo {author} {\bibfnamefont {T.}~\bibnamefont {Bartsch}},
  \bibinfo {author} {\bibfnamefont {F.}~\bibnamefont {Revuelta}}, \bibinfo
  {author} {\bibfnamefont {F.}~\bibnamefont {Borondo}}, \bibinfo {author}
  {\bibfnamefont {R.}~\bibnamefont {Benito}}, \ and\ \bibinfo {author}
  {\bibfnamefont {R.}~\bibnamefont {Hernandez}},\ }\href@noop {} {\bibfield
  {journal} {\bibinfo  {journal} {Physical Chemistry Chemical Physics}\
  }\textbf {\bibinfo {volume} {18}},\ \bibinfo {pages} {30270} (\bibinfo {year}
  {2016})}\BibitemShut {NoStop}%
\bibitem [{\citenamefont {Bardakcioglu}\ \emph {et~al.}(2018)\citenamefont
  {Bardakcioglu}, \citenamefont {Junginger}, \citenamefont {Feldmaier},
  \citenamefont {Main},\ and\ \citenamefont {Hernandez}}]{bardakcioglu2018}%
  \BibitemOpen
  \bibfield  {author} {\bibinfo {author} {\bibfnamefont {R.}~\bibnamefont
  {Bardakcioglu}}, \bibinfo {author} {\bibfnamefont {A.}~\bibnamefont
  {Junginger}}, \bibinfo {author} {\bibfnamefont {M.}~\bibnamefont
  {Feldmaier}}, \bibinfo {author} {\bibfnamefont {J.}~\bibnamefont {Main}}, \
  and\ \bibinfo {author} {\bibfnamefont {R.}~\bibnamefont {Hernandez}},\
  }\href@noop {} {\bibfield  {journal} {\bibinfo  {journal} {Physical Review
  E}\ }\textbf {\bibinfo {volume} {98}} (\bibinfo {year} {2018})}\BibitemShut
  {NoStop}%
\bibitem [{\citenamefont {Ezra}\ and\ \citenamefont
  {Wiggins}(2018)}]{ezra_2018}%
  \BibitemOpen
  \bibfield  {author} {\bibinfo {author} {\bibfnamefont {G.~S.}\ \bibnamefont
  {Ezra}}\ and\ \bibinfo {author} {\bibfnamefont {S.}~\bibnamefont {Wiggins}},\
  }\href@noop {} {\bibfield  {journal} {\bibinfo  {journal} {The Journal of
  Physical Chemistry A}\ }\textbf {\bibinfo {volume} {122}},\ \bibinfo {pages}
  {8354} (\bibinfo {year} {2018})}\BibitemShut {NoStop}%
\bibitem [{\citenamefont {Contopoulos}(1970)}]{contopoulos1970}%
  \BibitemOpen
  \bibfield  {author} {\bibinfo {author} {\bibfnamefont {G.}~\bibnamefont
  {Contopoulos}},\ }\href@noop {} {\bibfield  {journal} {\bibinfo  {journal}
  {The Astronomical Journal}\ }\textbf {\bibinfo {volume} {75}},\ \bibinfo
  {pages} {96} (\bibinfo {year} {1970})}\BibitemShut {NoStop}%
\bibitem [{\citenamefont {Heller}, \citenamefont {Stechel},\ and\ \citenamefont
  {Davis}(1980{\natexlab{b}})}]{heller1980}%
  \BibitemOpen
  \bibfield  {author} {\bibinfo {author} {\bibfnamefont {E.~J.}\ \bibnamefont
  {Heller}}, \bibinfo {author} {\bibfnamefont {E.~B.}\ \bibnamefont {Stechel}},
  \ and\ \bibinfo {author} {\bibfnamefont {M.~J.}\ \bibnamefont {Davis}},\
  }\href@noop {} {\bibfield  {journal} {\bibinfo  {journal} {The Journal of
  Chemical Physics}\ }\textbf {\bibinfo {volume} {73}},\ \bibinfo {pages}
  {4720} (\bibinfo {year} {1980}{\natexlab{b}})}\BibitemShut {NoStop}%
\bibitem [{\citenamefont {Davis}\ and\ \citenamefont
  {Heller}(1981)}]{davis1981}%
  \BibitemOpen
  \bibfield  {author} {\bibinfo {author} {\bibfnamefont {M.~J.}\ \bibnamefont
  {Davis}}\ and\ \bibinfo {author} {\bibfnamefont {E.~J.}\ \bibnamefont
  {Heller}},\ }\href@noop {} {\bibfield  {journal} {\bibinfo  {journal} {The
  Journal of Chemical Physics}\ }\textbf {\bibinfo {volume} {75}},\ \bibinfo
  {pages} {246} (\bibinfo {year} {1981})}\BibitemShut {NoStop}%
\bibitem [{\citenamefont {Martens}\ and\ \citenamefont
  {Ezra}(1987)}]{martens1987}%
  \BibitemOpen
  \bibfield  {author} {\bibinfo {author} {\bibfnamefont {C.~C.}\ \bibnamefont
  {Martens}}\ and\ \bibinfo {author} {\bibfnamefont {G.~S.}\ \bibnamefont
  {Ezra}},\ }\href@noop {} {\bibfield  {journal} {\bibinfo  {journal} {The
  Journal of Chemical Physics}\ }\textbf {\bibinfo {volume} {86}},\ \bibinfo
  {pages} {279} (\bibinfo {year} {1987})}\BibitemShut {NoStop}%
\bibitem [{\citenamefont {Wiggins}(2016)}]{wiggins_role_2016}%
  \BibitemOpen
  \bibfield  {author} {\bibinfo {author} {\bibfnamefont {S.}~\bibnamefont
  {Wiggins}},\ }\href@noop {} {\bibfield  {journal} {\bibinfo  {journal}
  {Regular and Chaotic Dynamics}\ }\textbf {\bibinfo {volume} {21}},\ \bibinfo
  {pages} {621} (\bibinfo {year} {2016})}\BibitemShut {NoStop}%
\bibitem [{\citenamefont {Contopoulos}\ \emph {et~al.}(1994)\citenamefont
  {Contopoulos}, \citenamefont {Farantos}, \citenamefont {Papadaki},\ and\
  \citenamefont {Polymilis}}]{contopoulos_1994}%
  \BibitemOpen
  \bibfield  {author} {\bibinfo {author} {\bibfnamefont {G.}~\bibnamefont
  {Contopoulos}}, \bibinfo {author} {\bibfnamefont {S.~C.}\ \bibnamefont
  {Farantos}}, \bibinfo {author} {\bibfnamefont {H.}~\bibnamefont {Papadaki}},
  \ and\ \bibinfo {author} {\bibfnamefont {C.}~\bibnamefont {Polymilis}},\
  }\href@noop {} {\bibfield  {journal} {\bibinfo  {journal} {Physical Review
  E}\ }\textbf {\bibinfo {volume} {50}},\ \bibinfo {pages} {4399} (\bibinfo
  {year} {1994})}\BibitemShut {NoStop}%
\bibitem [{\citenamefont {Farantos}(1998)}]{farantos_1998}%
  \BibitemOpen
  \bibfield  {author} {\bibinfo {author} {\bibfnamefont {S.~C.}\ \bibnamefont
  {Farantos}},\ }\href@noop {} {\bibfield  {journal} {\bibinfo  {journal}
  {Computer Physics Communications}\ }\textbf {\bibinfo {volume} {108}},\
  \bibinfo {pages} {240} (\bibinfo {year} {1998})}\BibitemShut {NoStop}%
\bibitem [{\citenamefont {Demian}\ and\ \citenamefont
  {Wiggins}(2017)}]{demian2017}%
  \BibitemOpen
  \bibfield  {author} {\bibinfo {author} {\bibfnamefont {A.~S.}\ \bibnamefont
  {Demian}}\ and\ \bibinfo {author} {\bibfnamefont {S.}~\bibnamefont
  {Wiggins}},\ }\href@noop {} {\bibfield  {journal} {\bibinfo  {journal}
  {International Journal of Bifurcation and Chaos}\ }\textbf {\bibinfo {volume}
  {27}},\ \bibinfo {pages} {1750225} (\bibinfo {year} {2017})}\BibitemShut
  {NoStop}%
\bibitem [{\citenamefont {Junginger}\ and\ \citenamefont
  {Hernandez}(2016{\natexlab{b}})}]{junginger2016uncovering}%
  \BibitemOpen
  \bibfield  {author} {\bibinfo {author} {\bibfnamefont {A.}~\bibnamefont
  {Junginger}}\ and\ \bibinfo {author} {\bibfnamefont {R.}~\bibnamefont
  {Hernandez}},\ }\href@noop {} {\bibfield  {journal} {\bibinfo  {journal} {The
  Journal of Physical Chemistry B}\ }\textbf {\bibinfo {volume} {120}},\
  \bibinfo {pages} {1720} (\bibinfo {year} {2016}{\natexlab{b}})}\BibitemShut
  {NoStop}%
\bibitem [{\citenamefont {Junginger}\ \emph
  {et~al.}(2017{\natexlab{a}})\citenamefont {Junginger}, \citenamefont
  {Duvenbeck}, \citenamefont {Feldmaier}, \citenamefont {Main}, \citenamefont
  {Wunner},\ and\ \citenamefont {Hernandez}}]{junginger2017chemical}%
  \BibitemOpen
  \bibfield  {author} {\bibinfo {author} {\bibfnamefont {A.}~\bibnamefont
  {Junginger}}, \bibinfo {author} {\bibfnamefont {L.}~\bibnamefont
  {Duvenbeck}}, \bibinfo {author} {\bibfnamefont {M.}~\bibnamefont
  {Feldmaier}}, \bibinfo {author} {\bibfnamefont {J.}~\bibnamefont {Main}},
  \bibinfo {author} {\bibfnamefont {G.}~\bibnamefont {Wunner}}, \ and\ \bibinfo
  {author} {\bibfnamefont {R.}~\bibnamefont {Hernandez}},\ }\href@noop {}
  {\bibfield  {journal} {\bibinfo  {journal} {The Journal of chemical physics}\
  }\textbf {\bibinfo {volume} {147}},\ \bibinfo {pages} {064101} (\bibinfo
  {year} {2017}{\natexlab{a}})}\BibitemShut {NoStop}%
\bibitem [{\citenamefont {Junginger}\ \emph
  {et~al.}(2017{\natexlab{b}})\citenamefont {Junginger}, \citenamefont {Main},
  \citenamefont {Wunner},\ and\ \citenamefont
  {Hernandez}}]{junginger2017variational}%
  \BibitemOpen
  \bibfield  {author} {\bibinfo {author} {\bibfnamefont {A.}~\bibnamefont
  {Junginger}}, \bibinfo {author} {\bibfnamefont {J.}~\bibnamefont {Main}},
  \bibinfo {author} {\bibfnamefont {G.}~\bibnamefont {Wunner}}, \ and\ \bibinfo
  {author} {\bibfnamefont {R.}~\bibnamefont {Hernandez}},\ }\href@noop {}
  {\bibfield  {journal} {\bibinfo  {journal} {Physical Review A}\ }\textbf
  {\bibinfo {volume} {95}} (\bibinfo {year} {2017}{\natexlab{b}})}\BibitemShut
  {NoStop}%
\bibitem [{\citenamefont {Lopesino}\ \emph {et~al.}(2015)\citenamefont
  {Lopesino}, \citenamefont {Balibrea}, \citenamefont {Wiggins},\ and\
  \citenamefont {Mancho}}]{lopesino_2015}%
  \BibitemOpen
  \bibfield  {author} {\bibinfo {author} {\bibfnamefont {C.}~\bibnamefont
  {Lopesino}}, \bibinfo {author} {\bibfnamefont {F.}~\bibnamefont {Balibrea}},
  \bibinfo {author} {\bibfnamefont {S.}~\bibnamefont {Wiggins}}, \ and\
  \bibinfo {author} {\bibfnamefont {A.~M.}\ \bibnamefont {Mancho}},\
  }\href@noop {} {\bibfield  {journal} {\bibinfo  {journal} {Communications in
  Nonlinear Science and Numerical Simulation}\ }\textbf {\bibinfo {volume}
  {27}},\ \bibinfo {pages} {40} (\bibinfo {year} {2015})}\BibitemShut {NoStop}%
\bibitem [{\citenamefont {Craven}\ and\ \citenamefont
  {Hernandez}(2016)}]{craven2016deconstructing}%
  \BibitemOpen
  \bibfield  {author} {\bibinfo {author} {\bibfnamefont {G.~T.}\ \bibnamefont
  {Craven}}\ and\ \bibinfo {author} {\bibfnamefont {R.}~\bibnamefont
  {Hernandez}},\ }\href@noop {} {\bibfield  {journal} {\bibinfo  {journal}
  {Physical Chemistry Chemical Physics}\ }\textbf {\bibinfo {volume} {18}},\
  \bibinfo {pages} {4008} (\bibinfo {year} {2016})}\BibitemShut {NoStop}%
\bibitem [{\citenamefont {Craven}\ and\ \citenamefont
  {Hernandez}(2015)}]{craven2015lagrangian}%
  \BibitemOpen
  \bibfield  {author} {\bibinfo {author} {\bibfnamefont {G.~T.}\ \bibnamefont
  {Craven}}\ and\ \bibinfo {author} {\bibfnamefont {R.}~\bibnamefont
  {Hernandez}},\ }\href@noop {} {\bibfield  {journal} {\bibinfo  {journal}
  {Physical review letters}\ }\textbf {\bibinfo {volume} {115}},\ \bibinfo
  {pages} {148301} (\bibinfo {year} {2015})}\BibitemShut {NoStop}%
\bibitem [{\citenamefont {Meiss}(1997)}]{meiss_average_1997}%
  \BibitemOpen
  \bibfield  {author} {\bibinfo {author} {\bibfnamefont {J.~D.}\ \bibnamefont
  {Meiss}},\ }\href@noop {} {\bibfield  {journal} {\bibinfo  {journal} {Chaos:
  An Interdisciplinary Journal of Nonlinear Science}\ }\textbf {\bibinfo
  {volume} {7}},\ \bibinfo {pages} {139} (\bibinfo {year} {1997})}\BibitemShut
  {NoStop}%
\bibitem [{\citenamefont {Garc{\'i}a-Garrido}\ \emph
  {et~al.}(2018)\citenamefont {Garc{\'i}a-Garrido}, \citenamefont {Curbelo},
  \citenamefont {Mancho}, \citenamefont {Wiggins},\ and\ \citenamefont
  {Mechoso}}]{garcia-garrido_2018}%
  \BibitemOpen
  \bibfield  {author} {\bibinfo {author} {\bibfnamefont {V.~J.}\ \bibnamefont
  {Garc{\'i}a-Garrido}}, \bibinfo {author} {\bibfnamefont {J.}~\bibnamefont
  {Curbelo}}, \bibinfo {author} {\bibfnamefont {A.~M.}\ \bibnamefont {Mancho}},
  \bibinfo {author} {\bibfnamefont {S.}~\bibnamefont {Wiggins}}, \ and\
  \bibinfo {author} {\bibfnamefont {C.~R.}\ \bibnamefont {Mechoso}},\
  }\href@noop {} {\bibfield  {journal} {\bibinfo  {journal} {Regular and
  Chaotic Dynamics}\ }\textbf {\bibinfo {volume} {23}},\ \bibinfo {pages} {551}
  (\bibinfo {year} {2018})}\BibitemShut {NoStop}%
\bibitem [{\citenamefont {{De Leon}}(1992)}]{DeLeon1992}%
  \BibitemOpen
  \bibfield  {author} {\bibinfo {author} {\bibfnamefont {N.}~\bibnamefont {{De
  Leon}}},\ }\href@noop {} {\bibfield  {journal} {\bibinfo  {journal} {J. Chem.
  Phys.}\ }\textbf {\bibinfo {volume} {96}},\ \bibinfo {pages} {285} (\bibinfo
  {year} {1992})}\BibitemShut {NoStop}%
\bibitem [{\citenamefont {{De Leon}}, \citenamefont {Mehta},\ and\
  \citenamefont {Topper}(1991{\natexlab{a}})}]{DeMeTo1991}%
  \BibitemOpen
  \bibfield  {author} {\bibinfo {author} {\bibfnamefont {N.}~\bibnamefont {{De
  Leon}}}, \bibinfo {author} {\bibfnamefont {M.~A.}\ \bibnamefont {Mehta}}, \
  and\ \bibinfo {author} {\bibfnamefont {R.~Q.}\ \bibnamefont {Topper}},\
  }\href@noop {} {\bibfield  {journal} {\bibinfo  {journal} {J. Chem. Phys.}\
  }\textbf {\bibinfo {volume} {94}},\ \bibinfo {pages} {8310} (\bibinfo {year}
  {1991}{\natexlab{a}})}\BibitemShut {NoStop}%
\bibitem [{\citenamefont {{De Leon}}, \citenamefont {Mehta},\ and\
  \citenamefont {Topper}(1991{\natexlab{b}})}]{DeMeTo1991a}%
  \BibitemOpen
  \bibfield  {author} {\bibinfo {author} {\bibfnamefont {N.}~\bibnamefont {{De
  Leon}}}, \bibinfo {author} {\bibfnamefont {M.~A.}\ \bibnamefont {Mehta}}, \
  and\ \bibinfo {author} {\bibfnamefont {R.~Q.}\ \bibnamefont {Topper}},\
  }\href@noop {} {\bibfield  {journal} {\bibinfo  {journal} {J. Chem. Phys.}\
  }\textbf {\bibinfo {volume} {94}},\ \bibinfo {pages} {8329} (\bibinfo {year}
  {1991}{\natexlab{b}})}\BibitemShut {NoStop}%
\bibitem [{\citenamefont {Feldmaier}\ \emph {et~al.}(2019)\citenamefont
  {Feldmaier}, \citenamefont {Schraft}, \citenamefont {Bardakcioglu},
  \citenamefont {Reiff}, \citenamefont {Lober}, \citenamefont {Tschöpe},
  \citenamefont {Junginger}, \citenamefont {Main}, \citenamefont {Bartsch},\
  and\ \citenamefont {Hernandez}}]{feldmaier_invariant_2019}%
  \BibitemOpen
  \bibfield  {author} {\bibinfo {author} {\bibfnamefont {M.}~\bibnamefont
  {Feldmaier}}, \bibinfo {author} {\bibfnamefont {P.}~\bibnamefont {Schraft}},
  \bibinfo {author} {\bibfnamefont {R.}~\bibnamefont {Bardakcioglu}}, \bibinfo
  {author} {\bibfnamefont {J.}~\bibnamefont {Reiff}}, \bibinfo {author}
  {\bibfnamefont {M.}~\bibnamefont {Lober}}, \bibinfo {author} {\bibfnamefont
  {M.}~\bibnamefont {Tschöpe}}, \bibinfo {author} {\bibfnamefont
  {A.}~\bibnamefont {Junginger}}, \bibinfo {author} {\bibfnamefont
  {J.}~\bibnamefont {Main}}, \bibinfo {author} {\bibfnamefont {T.}~\bibnamefont
  {Bartsch}}, \ and\ \bibinfo {author} {\bibfnamefont {R.}~\bibnamefont
  {Hernandez}},\ }\href@noop {} {\bibfield  {journal} {\bibinfo  {journal} {The
  Journal of Physical Chemistry B}\ }\textbf {\bibinfo {volume} {123}},\
  \bibinfo {pages} {2070} (\bibinfo {year} {2019})}\BibitemShut {NoStop}%
\bibitem [{Note2()}]{Note2}%
  \BibitemOpen
  \bibinfo {note} {{\protect \rm eq, top} and {\protect \rm eq, bot} in
  subscript denote the equilibrium points with positive y and negative y
  coordinates for the parametes chosen in this study}\BibitemShut {NoStop}%
\bibitem [{\citenamefont {Koon}\ \emph {et~al.}(2011)\citenamefont {Koon},
  \citenamefont {Lo}, \citenamefont {Marsden},\ and\ \citenamefont
  {Ross}}]{Koon2011}%
  \BibitemOpen
  \bibfield  {author} {\bibinfo {author} {\bibfnamefont {W.~S.}\ \bibnamefont
  {Koon}}, \bibinfo {author} {\bibfnamefont {M.~W.}\ \bibnamefont {Lo}},
  \bibinfo {author} {\bibfnamefont {J.~E.}\ \bibnamefont {Marsden}}, \ and\
  \bibinfo {author} {\bibfnamefont {S.~D.}\ \bibnamefont {Ross}},\ }\href@noop
  {} {\emph {\bibinfo {title} {{Dynamical systems, the three-body problem and
  space mission design}}}}\ (\bibinfo  {publisher} {Marsden books},\ \bibinfo
  {year} {2011})\ p.\ \bibinfo {pages} {327}\BibitemShut {NoStop}%
\bibitem [{\citenamefont {Parker}\ and\ \citenamefont
  {Chua}(1989)}]{Parker1989}%
  \BibitemOpen
  \bibfield  {author} {\bibinfo {author} {\bibfnamefont {T.~S.}\ \bibnamefont
  {Parker}}\ and\ \bibinfo {author} {\bibfnamefont {L.~O.}\ \bibnamefont
  {Chua}},\ }\href@noop {} {\emph {\bibinfo {title} {Practical Numerical
  Algorithms for Chaotic Systems}}}\ (\bibinfo  {publisher} {Springer-Verlag
  New York, Inc.},\ \bibinfo {address} {New York, NY, USA},\ \bibinfo {year}
  {1989})\BibitemShut {NoStop}%
\bibitem [{\citenamefont {Meyer}, \citenamefont {Hall},\ and\ \citenamefont
  {Offin}(2009)}]{Meyer2009}%
  \BibitemOpen
  \bibfield  {author} {\bibinfo {author} {\bibfnamefont {K.~R.}\ \bibnamefont
  {Meyer}}, \bibinfo {author} {\bibfnamefont {G.~R.}\ \bibnamefont {Hall}}, \
  and\ \bibinfo {author} {\bibfnamefont {D.}~\bibnamefont {Offin}},\
  }\href@noop {} {\emph {\bibinfo {title} {Applied Mathematical Sciences}}}\
  (\bibinfo  {publisher} {Springer},\ \bibinfo {year} {2009})\BibitemShut
  {NoStop}%
\bibitem [{\citenamefont {Marsden}\ and\ \citenamefont
  {Ross}(2006)}]{Marsden2006}%
  \BibitemOpen
  \bibfield  {author} {\bibinfo {author} {\bibfnamefont {J.~E.}\ \bibnamefont
  {Marsden}}\ and\ \bibinfo {author} {\bibfnamefont {S.~D.}\ \bibnamefont
  {Ross}},\ }\href@noop {} {\bibfield  {journal} {\bibinfo  {journal} {Bulletin
  of the American Mathematical Society}\ }\textbf {\bibinfo {volume} {43}},\
  \bibinfo {pages} {43} (\bibinfo {year} {2006})}\BibitemShut {NoStop}%
\end{thebibliography}%

\appendix

\section{Coupled harmonic 2 DoF system}
\label{ssect:coupled_2dof}

For this system, Hamilton's equations of motion are
\begin{equation}
\begin{aligned}
\dot{x} &= p_x \\
\dot{y} &= p_y \\
\dot{p}_x &= -\frac{\partial V}{\partial x} =  - (\omega_x^2 x + \delta y^2)  \\
\dot{p}_y &= -\frac{\partial V}{\partial y} =  - (\omega_y^2 y + 2\delta x y)
\end{aligned}
\label{eqn:two_dof_Barbanis}
\end{equation}

\textit{Hill's region and zero velocity curve \textemdash} 
The projection of energy surface into configuration space, $(x,y)$ plane, is the region of 
energetically possible motion for an energy $E(x,y,p_x,p_y) = e$. Let $M(e)$ denote this projection 
defined as
\begin{align}
M(e) = \left\{ (x,y) \, \vert \, V_B(x,y) \leqslant e \right\} 
\label{eqn:hills_region}
\end{align}
where $V_B(x,y)$ is the potential energy function~\eqref{eqn:Hamiltonian_Barbanis}. The projection~\eqref{eqn:hills_region} 
of energy surface is known in mechanics as the {\it Hill's region}. The boundary of $M(e)$ is known 
as the zero velocity curve, 
and plays an important role in placing bounds on the motion of a phase space point for a given total energy. The zero 
velocity curves are the locus of points in the $(x,y)$ plane where the kinetic energy, and hence the angular velocity vector 
vanishes, that is

\begin{align}
E(x,y,p_x,p_y) = e =& \frac{1}{2}\left( p_x^2 + p_y^2 \right) + V_B(x,y) \\
p_x^2 + p_y^2  =& 2(e - V_B(x,y)) = 0	
\end{align}

From Eqn.~\eqref{eqn:hills_region}, it is clear that the state is only able to move on the side of this curve for which the 
kinetic energy is positive. The other side of the curve, where the kinetic energy is negative and motion is impossible, is 
referred to as the energetically forbidden realm, and shown as gray region.

\textit{Symmetries of the equations of motion \textemdash} We note the symmetries in the 
system~\eqref{eqn:two_dof_Barbanis}, by substituting $(-y,-p_y)$ for 
$(y, p_y)$ which implies reflection about the 
$x-$axis and expressed as
\begin{equation}
s_y: (x,y,p_x,p_y,t) \rightarrow (x,-y,p_x,-p_y,t)
\end{equation}
Thus, if $(x(t),y(t),v_x(t),v_y(t))$ is a solution to~\eqref{eqn:two_dof_Barbanis}, then 
$(x(t),-y(t),p_x(t),-p_y(t))$ is another solution. The conservative system also has time-reversal 
symmetry 
\begin{equation}
s_t: (x,y,p_x,p_y,t) \rightarrow (x,y,-p_x,-p_y,-t)
\end{equation}
So, if $(x(t),y(t),p_x(t),p_y(t))$ is a solution to~\eqref{eqn:two_dof_Barbanis}, then 
$(x(-t),y(-t),-p_x(-t),-p_y(-t))$ is another solution. These symmetries can be used to decrease the 
number of computations, and to find special solutions. For example, any solution 
of~\eqref{eqn:two_dof_Barbanis} will evolve on the energy surface given 
by~\eqref{eqn:Hamiltonian_Barbanis}. For fixed energy, $E(x,y,p_x,p_y) = e$, there will be zero 
velocity curves corresponding to $V_B(x, y) = e$, the contours shown in 
Fig.~\ref{fig:pes_cont_Barbanis}. Any trajectory which touches the zero velocity curve at time 
$t_0$ must retrace its path in configuration space (i.e., $q = (x, y)$ 
space), 
\begin{equation}
q(-t + t_0) = q(t + t_0) \qquad \mathring{q}(-t + t_0) = -\mathring{q}(t + t_0)
\end{equation}

\subsection{Computing the NHIM and its invariant manifolds associated with the index-1 saddle}
\label{ssect:tube_mani} 
For the ease with discussing the geometry, we call the equilibrium with positive y-coordinate $\mathbf{x}_{\rm eq, top}$ and negative y-coordinate $\mathbf{x}_{\rm eq, bot}$.

\textbf{Select appropriate energy above the critical value~\textemdash~} For computation of manifolds that act as 
\emph{boundary} between the transition and non-transition trajectories, we select the total energy, $E$, above the critical 
value and so the excess energy 
$\Delta E > 0$. This excess energy can be arbitrarily large as long as the energy surface stays within the dynamical 
system's phase space bounds.

\textbf{Differential correction and numerical continuation for the NHIM~\textemdash~} 
We consider a procedure 
which computes periodic orbits around  in a relatively straightforward fashion. This procedure begins with small ``seed'' 
initial conditions obtained from the linearized equations of motion near $\mathbf{x}_{\rm eq, bot}$\footnote{{\rm eq, top} 
and {\rm eq, bot} in subscript denote the equilibrium points with positive y and negative y coordinates for the parametes 
chosen in this study} and uses differential correction and numerical continuation to generate the desired 
periodic orbit corresponding to the chosen energy $E$~(\cite{Koon2011}). The result is a periodic orbit of the desired 
energy $E$ of some period $T$, which will be close to $2\pi/\omega$ where $\pm i\omega$ is the imaginary pair of eigenvalues 
of the linearization around the saddle point. 

\textit{Guess initial condition of the periodic orbit \textemdash} The linearized equations of motion near an equilibrium  
point can be used to initialize a guess for the differential correction method. Let us select the equilibrium point, 
$\mathbf{x}_{\rm eq, bot}$. The linearization yields an eigenvalue problem $Av = \gamma v $, where $A$ is the Jacobian 
matrix~\eqref{eqn:jacobian_2dof} evaluated at the equilibrium point, $\gamma$ is the eigenvalue, and $v = [k_1, k_2, k_3, 
k_4]^T$ is the corresponding eigenvector. Thus, using the structure of $A$ from Eqn.~\eqref{eqn:jacobian_2dof} we can write 
\begin{equation}
\begin{aligned}
k_3 &= \gamma k_1 \\
k_4 &= \gamma k_2 \\
a k_1 + b k_2 &= \gamma k_3 \\
c k_1 + d k_2 &= \gamma k_4 
\end{aligned}
\end{equation}
where $a,b,c,d$ are entries in the Jacobian~\eqref{eqn:jacobian_2dof} and evaluated at the equilibrium point 
$(x_{\rm eq, bot},y_{\rm eq, bot}, 0, 0)$. So when $\gamma = \pm \lambda$, which correspond to the saddle directions of the 
equilibrium point, the corresponding 
eigenvectors are 
\begin{equation}
\begin{aligned}
u_1 &= [1,k_2,\lambda,\lambda k_2] \\
u_2 &= [1,k_2,-\lambda, -\lambda k_2] 
\end{aligned}
\end{equation}
and when $\gamma = \pm i \omega$, which correspond to the center 
directions, the corresponding eigenvectors are 
\begin{equation}
\begin{aligned}
w_1 &= [1,k_2, i\omega, i\omega k_2] \\
w_2 &= [1,k_2, -i\omega, -i\omega k_2]
\end{aligned}
\end{equation}
where, $k_2 = (\gamma^2 - a)/b$ is the constant 
depending on the eigenvalue, $\gamma$. Thus, the general solution of linearized equation of motion in 
Eqn.~\eqref{eqn:gen_sol_ind1saddle_2dof} can be used to initialize a guess for the periodic orbit for a small amplitude, 
$A_x \approx 10^{-4}$. The idea is to use the complex eigenvalue and the corresponding eigenvector to obtain a starting 
guess for the initial condition on the periodic orbit and its period $T_{\rm po}$, which should be close to $2\pi/\omega$ 
(generalization of Liapounov's theorem) and increase monotonically with excess energy, $\Delta E$. 

The initial condition for a periodic orbit of x-amplitude, $A_x > 0$ can be computed by letting $A_1 = A_2 = 0$ and $t = 0$ 
in Eqn.~\eqref{eqn:gen_sol_ind1saddle_2dof}, and $\beta = -A_x/2$ (this choice is made to get rid of factor 2) denotes a 
small amplitude in the general linear solution. 
Thus, using the eigenvector along the center direction we can guess the initial condition to be
\begin{equation}
\begin{aligned}
\bar{\mathbf{x}}_{\rm 0,g} =& \begin{pmatrix}
x_{\rm eq,bot},y_{\rm eq,bot},0,0
\end{pmatrix}^T + 
2Re(\beta w_1) \\
=& \begin{pmatrix}
x_{\rm eq,bot} - A_x, y_{\rm eq,bot} - A_x k_2, 0, 0
\end{pmatrix}^T
\end{aligned}
\end{equation}
%  

%\textbf{Differential correction.}
Without loss of generality, let us consider the bottom index-1 saddle equilibrium point (on the potential energy surface, $y 
< 0$), then the initial guess is given by   
\begin{equation}
\begin{aligned}
\bar{\mathbf{x}}_{0,\rm g} =& \left( -\frac{\omega_y^2}{2 \delta}, -\frac{1}{\sqrt{2}}\frac{\omega_x 
	\omega_y}{\delta}, 0, 0 \right) + \\
	& \left( A_x, \frac{A_x (\omega^2 - \omega_x^2)}{\sqrt{2}\omega_x 
	\omega_y}, 0, 0 \right)
\end{aligned}
\end{equation}

\textit{Differential correction of the initial condition \textemdash}~In this 
procedure, we attempt to introduce small change in the initial guess such that 
the periodic orbit $\bar{\mathbf{x}}_{\rm po}$
\begin{align}
\left\| \bar{\mathbf{x}}_{\rm po}(T) - 
\bar{\mathbf{x}}_{\rm po}(0) \right\| < \epsilon
\end{align}
for some tolerance $\epsilon << 1$. In this approach, we hold $x-$coordinate 
constant, while applying correction to the initial guess of the $y-$coordinate, 
use $v_y-$coordinate for terminating event-based integration, and
%for event of Poinca\`e map crossing during numerical integration and 
$v_x-$coordinate to test convergence of the periodic 
orbit. It is to be noted that this combination of coordinates is suitable for 
the structure of initial guess at hand, and in general will 
require some permutation of the phase space coordinates 
to achieve a stable algorithm. 

Let us denote the flow map of a differential equation 
$\mathring{\mathbf{x}} = \mathbf{f}(\mathbf{x})$ with 
initial condition $\mathbf{x}(t_0) = \mathbf{x}_0$ by 
$\phi(t;\mathbf{x}_0)$. Thus, the displacement 
of the final state under a perturbation $\delta t$ 
becomes 
\begin{align}
\delta \bar{\mathbf{x}}(t + \delta t) = \phi(t + \delta 
t;\bar{\mathbf{x}}_0 + \delta \bar{\mathbf{x}}_0) - 
\phi(t ;\bar{\mathbf{x}}_0)
\end{align}
with respect to the reference orbit 
$\bar{\mathbf{x}}(t)$. Thus, measuring the displacement 
at $t_1 + \delta t_1$ and expanding into Taylor series 
gives
\begin{align}
\delta \bar{\mathbf{x}}(t_1 + \delta t_1) = 
\frac{\partial \phi(t_1;\bar{\mathbf{x}}_0)}{\partial 
	\mathbf{x}_0}\delta \bar{\mathbf{x}}_0 + \frac{\partial 
	\phi(t_1;\bar{\mathbf{x}}_0)}{\partial t_1}\delta t_1 + 
h.o.t
\end{align}
where the first term on the right hand side is the state 
transition matrix, $\mathbf{\Phi}(t_1,t_0)$, when 
$\delta 
t_1 = 0$. Thus, it can be obtained as numerical solution 
to the variational equations as discussed in~\cite{Parker1989}. Let 
us suppose we want to reach the desired point 
$\mathbf{x}_{\rm d}$, we have
\begin{align}
\bar{\mathbf{x}}(t_1) = \phi(t_1;\bar{\mathbf{x}}_0) 
= \bar{\mathbf{x}}_1 = \mathbf{x}_d - \delta 
\bar{\mathbf{x}}_1
\end{align}
which has an error $\delta \bar{\mathbf{x}}_1$ and needs 
correction. This correction to the first order can be 
obtained from the state transition matrix at $t_1$ and an 
iterative procedure of this small correction based on first order yields 
convergence in few steps. For the equilibrium point under consideration, we 
initialize the guess as
\begin{align}
\bar{\mathbf{x}}(0) = (x_{0,g},y_{0,g},0,0)^T
\end{align}
and using numerical integrator we continue until next 
$v_x = 0$ event crossing with a high pecified tolerance 
(typically $10^{-14}$). So, we obtain 
$\bar{\mathbf{x}}(t_1)$ which for the guess periodic 
orbit denotes the half-period point, $t_1 = T_{0,g}/2$ 
and compute the state transition matrix 
$\mathbf{\Phi}(t_1,0)$. This can be used to correct the 
initial value of $y_{0,g}$ to approximate the periodic 
orbit while keeping $x_{0,g}$ constant. Thus, correction 
to the first order is given by
\begin{align}
\delta v_{x_1} = \Phi_{32}\delta y_0 + 
\mathring{v}_{x_1}\delta t_1 + h.o.t \\
\delta v_{y_1} = \Phi_{42}\delta y_0 + 
\mathring{v}_{y_1}\delta t_1 + h.o.t
\end{align}
where $\Phi_{ij}$ is the $(i,j)^{th}$ entry of 
$\mathbf{\Phi}(t_1,0)$ and the acceleration terms come 
from the equations of motion evaluated at the crossing 
$t 
= t_1$ when $v_{x_1} = \delta v_{x_1} = 0$. Thus, we 
obtain the first order correction $\delta y_0$ as
\begin{align}
\delta y_0 &\approx \left(\Phi_{42} - 
\Phi_{32}\frac{\mathring{v}_{y_1}}{\mathring{v}_{x_1}} 
\right)^{-1} \delta v_{y_1} \\
y_0 &\rightarrow y_0 - \delta y_0
\end{align}
which is iterated until $|v_{y_1}| = |\delta v_{y_1}| < 
\epsilon$ for some tolerance $\epsilon$, since we want 
the final periodic orbit to be of the form 
\begin{align}
\bar{\mathbf{x}}_{t_1} = (x_1,y_1,0,0)^T
\end{align}
This procedure yields an accurate initial condition for 
a periodic orbit of small amplitude $A_x << 1$, since 
our 
initial guess is based on the linear approximation near 
the equilibrium point. It is also to be noted that 
differential correction assumes the guess periodic orbit 
has a small error (for example in this system, of the order 
of $10^{-2}$) and can be corrected using first order form 
of the correction terms. If, however, larger steps in 
correction are applied this can lead to unstable convergence 
as the half-orbit overshoots between 
successive steps. Even though there are other algorithms for detecting unstable 
periodic orbits, differential correction is easy to implement and shows 
reliable convergence for generating family of periodic orbits at 
arbitrary high excess energy near the index-1 saddle. 
%These unstable perioidic orbit at high excess energy are required for 
%constructing codimension$-1$ invariant manifolds .

% 
\begin{figure}[!th]
	\centering
	\subfigure[]{\includegraphics[width=0.45\textwidth]{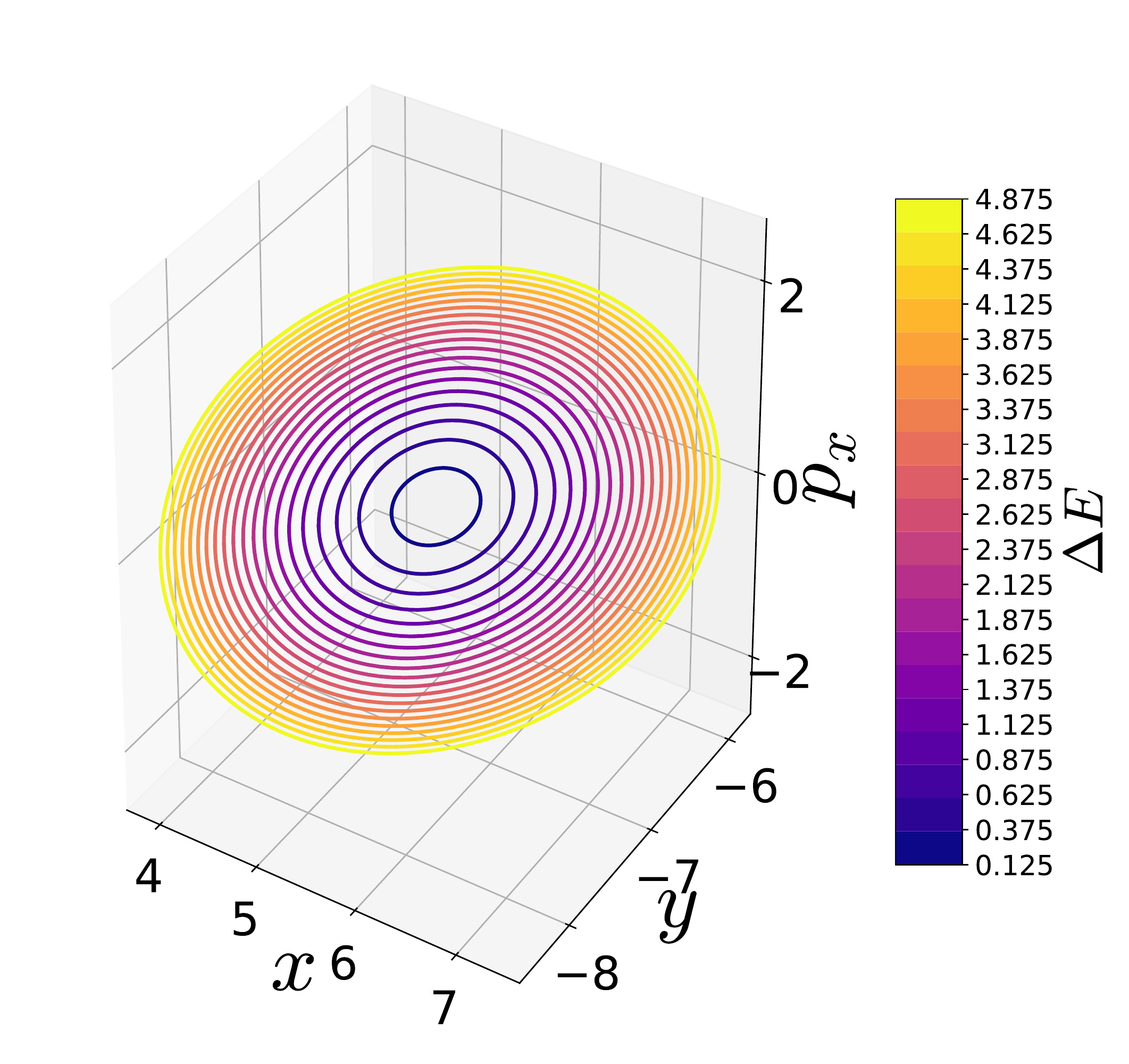}}
	\subfigure[]{\includegraphics[width=0.45\textwidth]{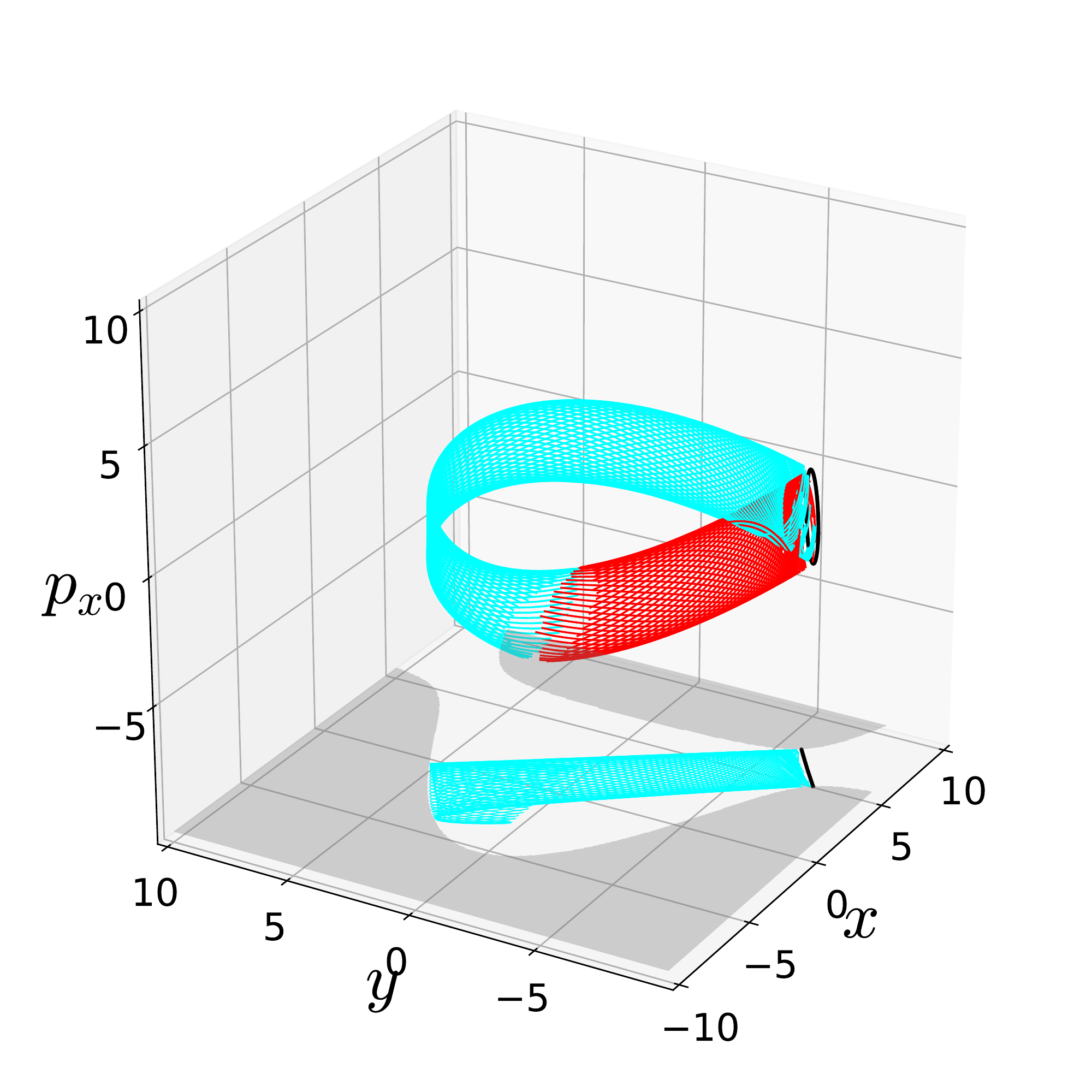}}
	\caption{Unstable periodic orbits at interval of $\Delta E = 0.25$ starting from $\Delta E = 0.125$ around the bottom 
	saddle equilibrium point. The stable and unstable manifolds associated with the unstable perioidic orbit (NHIM of 
	dimension 1) around the same equilibrium point at $\Delta E = 2.25$.}
	\label{fig:Barbanis_2dof_deltaE_upos_manifold}
\end{figure}

\textit{Numerical continuation to periodic orbit at arbitrary energy.\textemdash}~The procedure described above yields an 
accurate initial condition 
for a periodic orbit from a single initial guess. If our initial guess came 
from the linear approximation near the equilibrium point, from
Eqn.~\eqref{eqn:gen_sol_ind1saddle_2dof}, it has been observed numerically that we can only use this procedure for small 
amplitude, of order $10^{-4}$, periodic orbits 
around $\mathbf{x}_{\rm eq, bot}$. This small amplitude correspond to small excess 
energy, typically of the order $10^{-2}$, and if we want to compute the periodic orbit of arbitrarily large amplitude, we 
resort to numerical continuation for generating a family which reaches the appropriate total energy. This is done using two 
nearby periodic orbits of small amplitude to obtain initial guess for the next periodic orbit and performing differential 
correction to this guess. To this end, we proceed as follows. Suppose we find two small nearby periodic orbit initial 
conditions, $\bar{\mathbf{x}}_0^{(1)}$ and $\bar{\mathbf{x}}_0^{(2)}$, 
correct to within the tolerance $d_{\rm tol}$, using the differential correction 
procedure described above. We can generate a family of periodic orbits with 
successively increasing amplitudes around $\bar{\mathbf{x}}_{\rm eq, bot}$ in 
the following way. Let 
\begin{equation}
\Delta = \bar{\mathbf{x}}_0^{(2)} - \bar{\mathbf{x}}_0^{(1)} 
= [\Delta x_0, \Delta y_0, 0, 0]^T		
\end{equation}
A linear extrapolation to an initial guess of slightly larger amplitude, 
$\bar{\mathbf{x}}_0^{(3)}$ is given by
\begin{align}
\bar{\mathbf{x}}_{0,g}^{(3)} =&~\bar{\mathbf{x}}_0^{(2)} + \Delta \\
=& \left[(\mathbf{x}_0^{(2)} + \Delta x_0), (y_0^{(2)} + \Delta y_0), 0, 0  
\right] ^T \\
=& \left[x_0^{(3)}, y_0^{(3)}, 0, 0  \right] ^T
\end{align}
Thus, keeping $x_0^{(3)}$ fixed, we can use differential correction on this 
initial condition to compute an accurate solution $\bar{\mathbf{x}}_0^{(3)}$ 
from the initial guess $\bar{\mathbf{x}}_{\rm 0,g}^{(3)}$ and repeat the 
process until we have a family of solutions. We can keep track of the energy of 
each periodic orbit and when we have two solutions, $\bar{\mathbf{x}}_0^{\rm 
	(k)}$ and $\bar{\mathbf{x}}_0^{\rm (k+1)}$, 
whose energy brackets the appropriate energy, $E$, we can resort to combining 
bisection and differential correction to these two periodic orbits until we 
converge to the desired periodic orbit to within a specified 
tolerance. Thus, the result is a periodic orbit at desired total energy $E$ and of 
some period $T$ with an initial condition $X_0$. This is shown in Fig.~\ref{fig:Barbanis_2dof_deltaE_upos_manifold} for a 
series of excess energy at intervals of 0.250.

\textbf{Globalization of invariant manifolds \textemdash}~We find the local 
approximation to the 
unstable and stable manifolds of the periodic orbit from the eigenvectors of 
the monodromy matrix. 
Next, the local linear approximation of the unstable (or stable) manifold in 
the form of a state vector is integrated in the nonlinear equations of motion 
to produce the approximation of the unstable (or stable) manifolds. This 
procedure is known as \textit{globalization of the manifolds} and we proceed as 
follows.

First, the state transition matrix $\Phi(t)$ along the periodic orbit with initial condition $X_0$ can be obtained 
numerically by integrating the variational equations 
along with the equations of motion from $t = 0$ to $t = T$. This is known as the monodromy matrix $M = \Phi(T)$ and the 
eigenvalues can be computed numerically. 
For 
Hamiltonian systems (see~\cite{Meyer2009} for details), tells us that the four 
eigenvalues of $M$ are of the form
\begin{align}
\lambda_1 > 1, \qquad \lambda_2 = \frac{1}{\lambda_1}, \qquad \lambda_3 = 
\lambda_4 = 1
\end{align}
The eigenvector associated with eigenvalue $\lambda_1$ is in the unstable 
direction, 
the eigenvector associated with eigenvalue $\lambda_2$ is in the stable 
direction. Let 
$e^{s}(X_0)$ denote the normalized (to 1) stable eigenvector, and $e^{u}(X_0)$ 
denote 
the normalized unstable eigenvector. We can compute the manifold by 
initializing along 
these eigenvectors as:
\begin{equation}
X^s(X_0) = X_0 + \epsilon e^s(X_0)
\end{equation}
for the stable manifold at $X_0$ along the periodic orbit as
\begin{equation}
X^u(X_0) = X_0 + \epsilon e^u(X_0)
\end{equation}
for the unstable manifold at $X_0$. Here the small displacement from $X_0$ is 
denoted 
by $\epsilon$ and its magnitude 
should be small enough to be within the validity of the linear estimate, yet 
not so 
small that the time of flight becomes too large due to asymptotic nature of the 
stable 
and unstable manifolds. Ref.~\cite{Koon2011} suggests typical values of 
$\epsilon > 
0$ corresponding to nondimensional position displacements of magnitude around 
$10^{-6}$. 
By numerically integrating the unstable vector forwards in time, using both 
$\epsilon$ 
and $-\epsilon$, for the forward and backward branches respectively, we 
generate 
trajectories shadowing the two branches, $W^u_{+}$ and $W^u_{-}$, of the 
unstable 
manifold of the periodic orbit. Similarly, by integrating the stable vector 
backwards 
in time, using both $\epsilon$ and $-\epsilon$, for forward and backward branch 
respectively, we generate trajectories shadowing the stable manifold, 
$W^{s}_{+,-}$. 
For the manifold at $X(t)$, one can simply use the state transition matrix to 
transport 
the eigenvectors from $X_0$ to $X(t)$:
\begin{equation}
X^s(X(t)) = \Phi(t,0)X^s(X_0)
\end{equation}
It is to be noted that since the state transition matrix does not preserve the 
norm, the resulting vector must be normalized. The globalized invariant 
manifolds associated with index-1 saddles are known as 
Conley-McGehee tubes~(\cite{Marsden2006}). These tubes form the skeleton 
of transition dynamics by acting as conduits for the states inside them to 
travel between potential wells.
% 
%\begin{figure}[!ht]
%	\centering
%	\subfigure[]{\includegraphics[width=0.24\textwidth]{./energy_tube_U1neg_DelE10.eps}}
%	\subfigure[]{\includegraphics[width=0.24\textwidth]{./energy_tube_U1neg_DelE30.eps}}
%	\subfigure[]{\includegraphics[width=0.24\textwidth]{./energy_tube_U1neg_DelE50.eps}}
%	\subfigure[]{\includegraphics[width=0.24\textwidth]{./energy_tube_U1neg_DelE70.eps}}
%	\caption{Shows the intersection of tube (stable) 
%		manifold and energy surface with the Poincar\'e SOS, $U_1^-$, in black 
%		and magenta, respectively, for $\Delta E = 100, 300, 500, 700 \; 
%		{\rm (cm/s)^2}$. The 
%		trajectories that are in the cyan region for the given energy $\Delta 
%		E$ lead to imminent transition from 
%		the quadrant 1 to 2. A movie of increasing transition area on the 
%		Poincar\'e section, $U_1^-$, 
%		can be found \href{https://youtu.be/YZKYx0N9Zug}{here}.}
%	\label{fig:energy_tube_U1neg}
%\end{figure}
%

The computation of codimension-1 separatrix associated with the unstable periodic orbit around a index-1 saddle begins with 
the linearized equations of motion. This is obtained after a coordinate transformation to the saddle equilibrium point and 
Taylor expansion of the equations of motion. Keeping the first order terms in this expansion, we obtain the eigenvalues and 
eigenvectors of the linearized system. The eigenvectors corresponding to the center direction provide the starting guess for 
computing the unstable periodic orbits for small excess energy, $\Delta E << 1$, above the saddle's energy. This iterative 
procedure performs small correction to the starting guess based on the terminal condition of the periodic orbit until a 
desired tolerance is satisfied. This procedure is known as differential correction and generates unstable periodic orbits 
for small excess energy. Next, a numerical continuation is implemented to follow the small energy (amplitude) periodic 
orbits out to high excess energies. We apply this procedure to the Barbanis 2 DoF system in \S~\ref{sec:model_prob_2dof} to 
generate the unstable periodic orbit and its associated invariant manifolds as shown in 
Fig.~\ref{fig:Barbanis_2dof_deltaE_upos_manifold}.

% \subsubsection{Unstable periodic orbits.} 
\subsection{The Linearized Hamiltonian System} To find the linearized equations around the saddle equilibria with 
coordinates 
$(x_e,y_e,0,0)$, we need the quadratic terms of the Hamiltonian~\eqref{eqn:Hamiltonian_Barbanis} expanded about the 
equilibrium point. After making a coordinate change with $(x_e, y_e, 0, 0)$ as the origin, the quadratic terms of the 
Hamiltonian function for the linearized equations, which we shall call $H_l$, is given by
\begin{equation}
H_l = \frac{1}{2}p_x^2 + 
\frac{1}{2}p_y^2 + \frac{1}{2}\omega_x^2 x^2 + \frac{1}{2}\omega_y^2 y^2 +
2 \delta x y y_e + \delta y^2 x_e 
\label{eqn:Hamiltonian_2order}
\end{equation}

This gives the linear equations of motion near the equilibrium point as
\begin{equation}
\begin{pmatrix}
\dot{x} \\
\dot{y} \\
\dot{p_x}\\
\dot{p_y}
\end{pmatrix} =  \begin{pmatrix}
0 & 0 & 1 & 0 \\
0 & 0 & 0 & 1 \\
-\omega_x^2 & -2\delta y_e  & 0 & 0 \\
-2\delta y_e & -(\omega_y^2 + 2\delta x_e) & 0 & 0
\end{pmatrix}
\begin{pmatrix}
x \\
y \\
p_x\\
p_y
\end{pmatrix}
\end{equation}

\subsubsection{Linear analysis near the equilibria}
\label{ssect:linear}

Studying the linearization of the dynamics near the equilibria is an essential ingredient for understanding the more full 
nonlinear dynamics. We analyze the linearized dynamics near the saddle equilibrium points which extends to the full 
nonlinear system due to the generalization of Liapounov's theorem. Here we perform linearization of the vector 
field~\eqref{eqn:two_dof_Barbanis} to study the dynamics near the equilibrium points. This is given by the Jacobian, $D 
f(\mathbf{x})$, of the vector field
\begin{equation}
\mathbb{J} = D f(\mathbf{x}) =
\begin{pmatrix}
0 & 0 & 1 & 0 \\
0 & 0 & 0 & 1 \\
-\omega_x^2 & -2\delta y_e  & 0 & 0 \\
-2\delta y_e & -(\omega_y^2 + 2\delta x_e) & 0 & 0
\end{pmatrix}
\label{eqn:jacobian_2dof}
\end{equation}
where $\mathbf{x} = (x, y, p_x, p_y)$ and $(x_e, y_e, 0, 0)$ is the equilibrium point.

\textit{Saddle-center equilibrium. \textemdash} 
%\textbf{Index-1 saddle:} 
At the equilibria, 
\begin{equation}
\left(-\frac{\omega_y^2}{2\delta}, \pm 
\frac{1}{\sqrt{2}}\frac{\omega_x \omega_y}{\delta}, 0, 0 \right) 
\end{equation}
the Jacobian~\eqref{eqn:jacobian_2dof} becomes
\begin{equation}
\begin{aligned}
\mathbb{J}  =& D f(\mathbf{x})\Bigr|_{\left (-\frac{\omega_y^2}{2\delta}, \pm\frac{1}{\sqrt{2}}\frac{\omega_x 
		\omega_y}{\delta}, 0, 0 \right)} \\
=&
\begin{pmatrix}
0 & 0 & 1 & 0 \\
0 & 0 & 0 & 1 \\
-\omega_x^2 & \mp \sqrt{2}\omega_x \omega_y & 0 & 0 \\
\mp \sqrt{2}\omega_x \omega_y & 0 & 0 & 0
\end{pmatrix}	
\label{eqn:jacobian_saddles}
\end{aligned}
\end{equation}
The characteristic polynomial of the Jacobian~\eqref{eqn:jacobian_saddles} can be expressed as
\begin{equation}
{\rm det}(\mathbb{J} - \beta \mathbb{I}) = p(\beta) = \beta^4 + \omega_x^2 \beta^2 - 2\omega_x^2 \omega_y^2
\end{equation}
Let $\alpha = \beta^2$, then the roots of the above polynomial are $\sqrt{\alpha_1}, 
-\sqrt{\alpha_1}, \sqrt{\alpha_2}, -\sqrt{\alpha_2}$, where 

\begin{equation}
\begin{aligned}
\alpha_1 =& \frac{-\omega_x^2 + \sqrt{\omega_x^2 + 8\omega_x^2 \omega_y^2}}{2}, \\ 
\alpha_2 =& \frac{-\omega_x^2 - \sqrt{\omega_x^2 + 8\omega_x^2 \omega_y^2}}{2} 
\label{eqn:eigenvals_jac}
\end{aligned}
\end{equation}    
It is clear that $\alpha_1 > 0$ and $\alpha_2 < 0$, so let us define, $\lambda = \sqrt{\alpha_1}$ 
and $\omega = \sqrt{-\alpha_2}$, so the eigenvalues are $\lambda, -\lambda, i 
\omega, - i \omega$. This implies the equilibrium point is of type saddle $\times$ center type, or 
{\it index-1 saddle}.

%As shown in \S~\ref{subsect:linear}, 
As shown above, the eigenvalues are of the form: $\lambda, -\lambda, i \omega, - i \omega$ where
\begin{equation}
\begin{aligned}
\lambda^2 = \alpha_1 =& \frac{-\omega_x^2 + \sqrt{\omega_x^2 + 8\omega_x^2 \omega_y^2}}{2} \\
\omega^2 =  -\alpha_2 =& \frac{\omega_x^2 + \sqrt{\omega_x^2 + 8\omega_x^2 \omega_y^2}}{2} 
\end{aligned}
\end{equation}
and the eigenvectors as obtained in~\S~\ref{sect:eigvec_lin_saddle} are
\begin{equation}
\begin{aligned}
u_{\lambda} = \begin{bmatrix}
1, & -\frac{\lambda^2 + \omega_x^2}{\sqrt{2}\omega_x \omega_y}, & \lambda, & - 
\frac{\lambda(\lambda^2 + \omega_x^2)}{\sqrt{2}\omega_x \omega_y} 
\end{bmatrix} \\
u_{-\lambda} = \begin{bmatrix}
1, & -\frac{\lambda^2 + \omega_x^2}{\sqrt{2}\omega_x \omega_y}, & -\lambda, &  
\frac{\lambda(\lambda^2 + \omega_x^2)}{\sqrt{2}\omega_x \omega_y} 
\end{bmatrix}
\end{aligned}
\end{equation}
and 
\begin{equation}
\begin{aligned}
w_{i \omega} =& \begin{bmatrix}
1, & \frac{\omega^2 - \omega_x^2}{\sqrt{2}\omega_x \omega_y}, & i \omega, & \frac{i \omega( 
	\omega^2 - \omega_x^2)}{\sqrt{2}\omega_x \omega_y}
\end{bmatrix}\\
w_{-i \omega} =& \begin{bmatrix}
1, & \frac{\omega^2 - \omega_x^2}{\sqrt{2}\omega_x \omega_y}, & -i\omega, & 
-\frac{i \omega(\omega^2 - \omega_x^2)}{\sqrt{2}\omega_x \omega_y}
\end{bmatrix}
\end{aligned}
\end{equation}

\subsubsection{Eigenvectors of linearized system near index-1 saddle}
\label{sect:eigvec_lin_saddle}

Let us assume the eigenvector to be $v = [k_1, k_2, k_3, k_4]^T$, and the eigenvalue problem 
becomes $\mathbb{J}v = \beta v$. This gives the expressions
\begin{align}
k_3 =& \beta k_1 \label{eqn:row1} \\
k_4 =& \beta k_2 \label{eqn:row2} \\
-\omega_x^2 k_1 - \sqrt{2}\omega_x \omega_y k_2 =& \beta k_3 \label{eqn:row3}\\
-\sqrt{2} \omega_x \omega_y k_1 =& \beta k_4 \label{eqn:row4}
\end{align}
Let $k_1 = 1$, then using Eqns.~\eqref{eqn:row1}and~\eqref{eqn:row2} the eigenvector becomes 
$\bigl[1, k_2, \beta, \beta k_2 \bigr]$. 

Then, using Eqns.~\eqref{eqn:row3}and~\eqref{eqn:row4} for 
eigenvalue $\beta = \pm \lambda$, we get

\begin{equation}
\begin{aligned}
-\omega_x^2 - \sqrt{2}\omega_x \omega_y k_2 =& \lambda^2 \\
-\sqrt{2} \omega_x \omega_y =& \lambda k_2  \\
-\omega_x^2 - \sqrt{2}\omega_x \omega_y k_2^\prime =& \lambda^2 \\
-\sqrt{2} \omega_x \omega_y =& \lambda k_2^\prime
\end{aligned}    
\end{equation}

These imply, $k_2 = k_2^\prime$ and 
\begin{equation}
k_2 = -\frac{\lambda^2 + \omega_x^2}{\sqrt{2}\omega_x \omega_y}
\end{equation}
A similar approach for the eigenvalues, $\beta = \pm \omega$, gives us
\begin{equation}
k_2 = \frac{\omega^2 - \omega_x^2}{\sqrt{2} \omega_x \omega_y}
\end{equation}
Thus, the eigenvectors corresponding to $\pm \lambda$ are
\begin{equation}
\begin{aligned}
u_{\lambda} = \begin{bmatrix}
1, & -\frac{\lambda^2 + \omega_x^2}{\sqrt{2}\omega_x \omega_y}, & \lambda, & - 
\frac{\lambda(\lambda^2 + \omega_x^2)}{\sqrt{2}\omega_x \omega_y} 
\end{bmatrix} \\
u_{-\lambda} = \begin{bmatrix}
1, & -\frac{\lambda^2 + \omega_x^2}{\sqrt{2}\omega_x \omega_y}, & -\lambda, &  
\frac{\lambda(\lambda^2 + \omega_x^2)}{\sqrt{2}\omega_x \omega_y} 
\end{bmatrix}
\end{aligned}
\end{equation}
and for the eigenvalues, $\pm i \omega$, we obtain
\begin{equation}
\begin{aligned}
w_{i \omega} =& \begin{bmatrix}
1, & \frac{\omega^2 - \omega_x^2}{\sqrt{2}\omega_x \omega_y}, & i \omega, & \frac{i \omega( 
	\omega^2 - \omega_x^2)}{\sqrt{2}\omega_x \omega_y}
\end{bmatrix}\\
w_{-i \omega} =& \begin{bmatrix}
1, & \frac{\omega^2 - \omega_x^2}{\sqrt{2}\omega_x \omega_y}, & -i\omega, & 
-\frac{i \omega(\omega^2 - \omega_x^2)}{\sqrt{2}\omega_x \omega_y}
\end{bmatrix}
\end{aligned}
\end{equation}
where $\lambda$ and $\omega$ are positive constants~\eqref{eqn:eigenvals_jac} that depend on the 
parameters of the potential energy surface. Thus, the general solution of the linear system near 
the saddle equilibrium point is given by
% 
% \begin{widetext}
\begin{equation}
\begin{aligned}
\mathbf{x}(t) = & \left\{  x(t), y(t), v_x(t), v_y(t)  \right\} \\
= & A_1 e^{\lambda t} u_{\lambda} + A_2 
e^{-\lambda t} u_{-\lambda} + 2{\rm Re}\left( \beta e^{i \omega t} w_{i \omega} \right)
\label{eqn:gen_sol_ind1saddle_2dof}
\end{aligned}
\end{equation}
% \end{widetext}

%
with $A_1, A_2$ being real and $\beta = \beta_1 + i \beta_2$ being complex.

\section{Coupled harmonic 3 DoF system}
\label{sect:coupled_3dof}

For this system, Hamilton's equations of motion are given by
%The equations of motion in original coordinates are given by
% 
% \begin{widetext}
\begin{equation}
\begin{aligned}
\dot{x} &= \frac{\partial V_{\rm BC}}{\partial p_x} = p_x \\
\dot{y} &= \frac{\partial V_{\rm BC}}{\partial p_y} = p_y \\
\dot{z} &= \frac{\partial V_{\rm BC}}{\partial p_z} = p_z \\
\dot{p}_x &= -\frac{\partial V_{\rm BC}}{\partial x} =  - (\omega_x^2 x - 2 \epsilon x y - 2 \eta x z)  \\
\dot{p}_y &= -\frac{\partial V_{\rm BC}}{\partial y} =  - (\omega_y^2 y - \epsilon x^2) \\
\dot{p}_z &= -\frac{\partial V_{\rm BC}}{\partial z} =  - (\omega_z^2 z - \eta x^2) 
\end{aligned}
\label{eqn:three_dof_Barbanis}
\end{equation}
% \end{widetext}
%

\textit{Linear analyis near equilibrium point \textemdash}~Here we study the dynamics near the equilibrium points using 
linearization of the vector field~\eqref{eqn:three_dof_Barbanis} given by the Jacobian, $D f(\mathbf{x})$, of the vector 
field evaluated at the equilibrium point, $\mathbf{x}_{\rm eq} = (x_e, y_e, z_e, 0, 0, 0)$.
\begin{equation}
\begin{aligned}
\dot{\mathbf{x}} = &  \mathbb{J}\Bigr|_{\left (x_e, y_e, z_e, 0, 0, 0 \right)} \mathbf{x} = D f(\mathbf{x})\Bigr|_{\left 
(x_e, y_e, z_e, 0, 0, 0 \right)} \mathbf{x}  \\ 
\dot{\mathbf{x}} = &
\begin{pmatrix}
0 & 0 & 0 & 1 & 0 & 0\\
0 & 0 & 0 & 0 & 1 & 0\\
0 & 0 & 0 & 0 & 0 & 1\\
-(\omega_x^2 - 2 \epsilon y_e - 2 \eta z_e) & 2\epsilon x_e  & 2 \eta x_e & 0 & 0 & 0 \\
2 \epsilon x_e & - \omega_y^2 & 0 & 0 & 0 & 0 \\
-2 \eta x_e & 0 & -\omega_z^2 & 0 & 0 & 0 \\
\end{pmatrix}
\begin{bmatrix}
x \\
y \\
z \\
p_x \\
p_y \\
p_z
\end{bmatrix}
\label{eqn:jacobian_Barbanis3dof}
\end{aligned}
\end{equation}

%
%The Jacobian~\eqref{eqn:jacobian_Barbanis3dof} at the equilibrium point $\left( 0, 0, 0, 0, 0, 0 \right)$ becomes
%%
%\begin{equation}
%\mathbb{J}\Bigr|_{\left (0, 0, 0, 0, 0, 0 \right)} = D f(\mathbf{x})\Bigr|_{(0,0,0,0,0,0)} =
%\begin{pmatrix}
%0 & 0 & 0 & 1 & 0 & 0\\
%0 & 0 & 0 & 0 & 1 & 0\\
%0 & 0 & 0 & 0 & 0 & 1\\
%-\omega_x^2 & 0  & 0 & 0 & 0 & 0\\
%0 & -\omega_y^2 & 0 & 0 & 0 & 0 \\
%0 & 0 & -\omega_z^2 & 0 & 0 & 0 \\
%\end{pmatrix}
%\label{eqn:jacobian_Barbanis3dof_bot}
%\end{equation}
%

\textit{Symmetries of the equations of motion~\textemdash~}We note the symmetries in the 
system~\eqref{eqn:three_dof_Barbanis}, by substituting $(-x,-p_x)$ for $(x, p_x)$ which implies reflection about the $x = 0$ 
plane and expressed as
\begin{equation}
s_x: (x,y,z,p_x,p_y,p_z,t) \rightarrow (-x,y,z,-p_x,p_y,p_z,t)
\end{equation}
Thus, if $(x(t),y(t),z(t),p_x(t),p_y(t),p_z(t))$ is a solution to~\eqref{eqn:three_dof_Barbanis}, then \\ 
$(-x(t),y(t),-p_x(t),p_y(t),p_z(t))$ is another solution. The conservative system also has time-reversal 
symmetry 
\begin{align}
s_t: (x,y,z,p_x,p_y,p_z,t) & \rightarrow (x,y,z,-p_x,-p_y,-p_z,t)
\end{align}
So, if $(x(t),y(t),z(t),p_x(t),p_y(t),p_z(t))$ is a solution to~\eqref{eqn:three_dof_Barbanis}, then \\
$(x(-t),y(-t),z(-t),-p_x(-t),-p_y(-t),-p_z(-t))$ is another solution. These symmetries will be used to decrease the 
number of computations, and to find special solutions. For example, any solution of~\eqref{eqn:three_dof_Barbanis} will evolve 
on the energy surface given by~\eqref{eqn:Hamiltonian_BC_3dof}. For fixed energy, $E(x,y,z,p_x,p_y,p_z) = e$, there will be 
zero velocity curves corresponding to $V_{BC}(x, y, z) = e$, the contours shown in Fig.~\ref{fig:pes_cont_Barbanis}. Any 
trajectory which touches the zero velocity curve at time $t_0$ must retrace its path in configuration space (i.e., $q = (x, 
y, z)$ space), 
\begin{equation}
q(-t + t_0) = q(t + t_0) \qquad \mathring{q}(-t + t_0) = -\mathring{q}(t + t_0)
\end{equation}

\end{document}